\documentclass{amsart}

\usepackage{latexsym,amsthm,amssymb}

\usepackage[leqno]{amsmath}

\usepackage{mathrsfs}

\usepackage[all]{xy}

\newtheorem{theorem}{Theorem}[section]
\newtheorem{lemma}[theorem]{Lemma}

\newtheorem{proposition}[theorem]{Proposition}
\newtheorem{cor}[theorem]{Corollary}
\newtheorem{remark}[theorem]{Remark}
\newtheorem{noname}[theorem]{}

\numberwithin{equation}{theorem}

\newenvironment{Proof}
{\noindent \emph{Proof.}}
{\hfill $\square$ \medskip}

\renewcommand{\mathcal}{\mathscr}

\newcommand{\SA}{{\mathcal{A}}}
\newcommand{\SB}{{\mathcal{B}}}
\newcommand{\SC}{{\mathcal{C}}}
\newcommand{\SD}{{\mathcal{D}}}
\newcommand{\SE}{{\mathcal{E}}}
\newcommand{\SF}{{\mathcal{F}}}
\newcommand{\SG}{{\mathcal{G}}}
\newcommand{\SH}{{\mathcal{H}}}
\newcommand{\SI}{{\mathcal{I}}}
\newcommand{\SJ}{{\mathcal{J}}}
\newcommand{\SK}{{\mathcal{K}}}
\newcommand{\SL}{{\mathcal{L}}}
\newcommand{\SM}{{\mathcal{M}}}
\newcommand{\SN}{{\mathcal{N}}}
\newcommand{\SO}{{\mathcal{O}}}
\newcommand{\SP}{{\mathcal{P}}}

\newcommand{\SU}{{\mathcal{U}}}
\newcommand{\SV}{{\mathcal{V}}}

\newcommand{\SX}{{\mathcal{X}}}
\newcommand{\SY}{{\mathcal{Y}}}

\newcommand{\PP}{\mathbb{P}}

\newcommand{\WY}{\widetilde{Y}}
\newcommand{\WX}{\widetilde{X}}

\newcommand{\WJ}{\widetilde{J}}
\newcommand{\WL}{\widetilde{L}}

\newcommand{\wl}{\widetilde{l}}

\newcommand{\wi}{\widetilde{i}}
\newcommand{\wpp}{\widetilde{p}}
\newcommand{\wphi}{\widetilde{\varphi}}
\newcommand{\whphi}{\widehat{\varphi}}

\newcommand{\womega}{\widetilde{\omega}}

\newcommand{\smV}{\scriptscriptstyle{V}}

\newcommand{\smU}{\scriptscriptstyle{U}}

\title{Smoothing of ribbons over curves}

\author{Miguel Gonz\'alez}

\address{Departamento de \'Algebra, Facultad de Ciencias Matem\'aticas,
Universidad Complutense de Madrid.}

\email{mgonza@mat.ucm.es}

\thanks{{\em Acknowledgement.}
I would like to express my heartfelt thanks to my Ph.D. advisor Francisco~Javier Gallego for many hints and helpful conversations, for his patient encouragement and generous assistance along the work of this paper.\\
This paper was supported by the research project BFM2000--0621 from Spanish Ministry of Science and Technology.}
\begin{document}
\begin{abstract}
We prove that ribbons, i.e. double structures associated with a line bundle $\SE$ over its reduced support, a smooth irreducible projective curve of arbitrary genus, are smoothable if their arithmetic genus is greater than or equal to $3\,$ and the support curve possesses a smooth irreducible double cover with trace zero module $\SE$.
The method we use is based on the infinitesimal techniques that we develop to show that if the support curve admits such a double cover then every embedded ribbon over the curve is ``infinitesimally smoothable'', i.e. the ribbon can be obtained as central fiber of the image of some first--order infinitesimal deformation of the map obtained by composing the double cover with the embedding of the reduced support in the ambient projective space containing the ribbon.
We also obtain embeddings in the same projective space for all ribbons associated with $\SE$.
Then, assuming the existence of the double cover, we prove that the ``infinitesimal smoothing'' can be extended to a global embedded smoothing for embedded ribbons of arithmetic genus greater than or equal to $3$.
As a consequence we obtain the smoothing results.
\end{abstract}
\maketitle
\section*{Introduction}
A ribbon $\WY$ is a multiplicity two structure associated with a line bundle $\SE$ over its reduced support~$Y$.
Precisely, $\SE$ is the ideal $\SI$ of $Y$ inside $\WY$, (from $\SI^2=0$ it follows that $\SI$ is an $\SO_Y$--module).
The scheme $\WY$ is called a ribbon over $Y$ with conormal bundle $\SE$,~\cite[\S1]{BayerEisenbud95}.
The notion of ribbon can be extended to higher multiplicity by allowing $\SE$ to have any rank.
If $\SE$, instead of having rank $1$, has rank $n-1$ then $\WY$ is called a rope of multiplicity $n$.\\
A smoothing of a ribbon is a family, flat over a smooth pointed affine curve, whose general fiber is a smooth variety an whose special fiber is the ribbon.
If the ribbon is embedded in an ambient variety and the family is a subvariety of the product of the ambient variety and the base curve of the family, then we call it an embedded smoothing.\\
Ribbons were first studied at length by D.~Bayer and D.~Eisenbud in their fundamental work~\cite{BayerEisenbud95}.
They are important as far as they appear as degenerations of smooth varieties whose properties are interesting to study.
Often, those properties are easier to study on a degenerate variety which, nevertheless, has a simpler structure in many ways (e.g. the structure of the Picard group of a ``rational" ribbon, i.e. a ribbon over $\PP^1$, or the computation of the equations of a $K3$--carpet, i.e. a ribbon with reduced support on a smooth rational normal scroll and conormal bundle the canonical line bundle of the scroll, see~\cite{BayerEisenbud95}), though its nondegenerate counterpart, which even having nicer properties from a geometric point of view (smoothness, irreducibility, ...), has a more complex structure.
Indeed, one of the source of interest for ribbons and other double structures in the $90$'s was the study of Green's conjecture regarding the syzygies of a canonical curve.
For this approach to be effective one needs to know, to start with, that ribbons are smoothable.
That is the reason why finding ways of smoothing a ribbon is important.\\
The goal of this paper is to determine under what conditions a ribbon can be smoothed. The result we obtain is that over a smooth irreducible projective curve $Y$ of arbitrary genus $g$, every ribbon (with very few exceptions if $g=0$ or $g=1$) with conormal bundle $\SE$ is smoothable under some natural geometric condition. This geometric condition is that there is a double cover of $Y$ with trace zero module $\SE$.\\

The appearance of ribbons as flat limit of smooth curves is expected whenever a family of embeddings degenerates into a $2:1$ morphism over a smooth curve $Y$. Indeed, in this situation the degree and the genus imposed by the images of the embeddings on their flat limit indicate that this flat limit must be a double structure over $Y$, whose genus is the genus of a ribbon over $Y$ associated with the trace zero module of the double cover.
This situation leads us to think that the natural geometric condition that we need to impose if we expect ribbons to be smoothable is that there is a double cover of $Y$ with trace zero module $\SE$.\\
At the infinitesimal level the fact that an embedded ribbon is contained in the first infinitesimal neighborhood of its reduced support inside the ambient variety suggests that the ribbon can be captured by ``first--order infinitesimal smoothings".
More specifically, one would like to establish the following correspondence: on one side we would have an embedded ribbon $\WY$, on the other we would have a first--order infinitesimal deformation of the composite map obtained from the double cover of $Y$ and the embedding of $Y$ in the ambient variety, in such a way that the central fiber of the image of this deformation is the ribbon.
This correspondence does not come totally unexpected (see~\cite{Fong93}, where canonically embedded non--hyperelliptic rational ri\-bbons of arithmetic genus greater than or equal to $3$ appear associated with infinitesimal deformations of hyperelliptic covers of rational normal curves, and~\cite{GallegoPurna97}, where $K3$--carpets appear associated with deformations of hyperelliptic covers of rational normal scrolls).
These examples are, nevertheless, particular cases and the approach is in one case done by an explicit computation (an approach that one hope to work only when the reduced curve $Y$ is as simple as $\PP^1$ and the line bundle associated with the map from the double cover to the ambient projective space is well--behaved as is the case with the canonical line bundle).
In the case of $K3$--carpets, the approach is not an explicit computation of the equations of an infinitesimal deformation of a double cover, but the success of the proof relies heavily on some very special characteristics, such as the existence of a unique double structure with $K3$ invariants on a given rational normal scroll.
Thus a more general and conceptual approach is needed.
This is what we do in this paper, setting up the foundations for the process on how a first--order infinitesimal deformation of the composite morphism obtained from a double cover of a curve $Y$ of arbitrary genus and the inclusion of $Y$ in the ambient projective space, produces a ribbon $\WY$ on $Y$, and how every ribbon $\WY$ on $Y$ comes indeed from such a process.
This is done in Proposition~\ref{construction} and Theorem~\ref{image} which say that every first--order, locally trivial, infinitesimal deformation of a morphism $X \overset{\varphi}\to Z$, which is finite over its image $Y$, produces a rope over $Y$ mapping to the target variety $Z$ and a first--order deformation of $Y$ in $Z$ so that the central fiber of the image of the deformation morphism is equal to the image of the rope and the whole image of the deformation morphism is the scheme-theoretic union of its central fiber and the flat embedded deformation.
These results give geometric content to the arrow, obtained by cohomological methods, from $H^0(\SN_{\varphi})$, the space of first--order, locally trivial, infinitesimal deformation of $\varphi$, to $\mathrm{Hom}(\mathcal{I}_{Y,Z}/\mathcal{I}_{Y,Z}^2,\SO_Y)\oplus \mathrm{Hom}(\mathcal{I}_{Y,Z}/\mathcal{I}_{Y,Z}^2,\SE)$.
The conceptual understanding of this process is also done in Theorem~\ref{main} (which can be understood as an ``infinitesimal smoothing'' result for ribbons) which says that if this arrow is surjective then every rope, with conormal bundle $\SE=\pi_*\SO_X/\SO_Y$ (where $\pi$ is the map from $X$ to $Y$), embedded in $Z\,$ is the central fiber of the image of some first--order, locally trivial, infinitesimal deformation of $\varphi$, and in Theorem~\ref{key-curves} which says that in the case where $X$ is a curve every rope over the curve $Y$ with conormal bundle $\SE$ is obtained via this process.
The way we prove the smoothing of ribbons is by showing that we can obtain the ribbon as central fiber of the image of a deformation of a morphism $2:1$ to a family of embeddings.
Our main smoothing result is Theorem~\ref{abssmoothing} which says:
\begin{theorem}\label{abssmoothing.0}
Let $Y$ be a smooth irreducible projective curve and let $\SE$ be a line bundle on $Y$.
Assume that there is a smooth irreducible double cover $X \overset{\pi} \to Y$ with $\pi_*\SO_X/\SO_Y = \SE$.
Then every ribbon $\WY$ over $Y$ with conormal bundle $\SE$ and arithmetic genus $p_a(\WY) \geq 3$ is smoothable.
\end{theorem}
As we said above, the existence of a double cover with trace zero module $\SE$ is the natural condition we might impose for the ribbons with conormal bundle $\SE$ to be smoothable.
In  fact, this condition turns out to be hardly restrictive in comparison with the obvious necessary condition for a ribbon to be smoothable, namely: its arithmetic genus is greater than or equal to zero.
Besides, since the arithmetic genus of a ribbon $\WY$ with conormal bundle $\SE$ is $p_a(\WY)= d+2g-1$, where $d = - \mathrm{deg} \,\SE$ and $g$ is the genus of $Y$, the existence of such a double cover implies the condition $p_a(\WY) \geq 3$, with very few exceptions if $g=0$ or $g=1$.\\   
To get an idea of the scope of our results we point out that if the genus $g$ of $Y$ or the arithmetic genus $p_a$ of $\WY$ gets bigger there exist more ribbons of arithmetic genus $p_a$, over a curve $Y$ of genus $g$.
The reason is that the ribbons with conormal bundle $\SE$ over a curve $Y$ are classified, up to isomorphism over $Y$, by the elements of the space $\mathrm{Ext}_Y^1(\omega_Y, \SE)$, up to the action of $k^*$ (see~~\cite[1.2]{BayerEisenbud95} or~\ref{noname1.2}).
Notice that if $p_a \geq 2$, then the space $\mathrm{Ext}_Y^1(\omega_Y, \SE)$ has dimension $g-2+p_a$.\\\\
The smoothing of ribbons is obtained as an embedded smoothing since the smoothing appears as the image of a deformation of a morphism composite of a double cover of $Y$ and the embedding of $Y$ in a projective space.
This is done in Theorem~\ref{embsmoothing} which says:
\begin{theorem}\label{embsmoothing.0}
Let $Y$ be a smooth irreducible projective curve and let $\SO_Y(1)$ be a very ample line bundle on $Y$.
Let $\SE$ be a line bundle on $Y$.
Assume that $\SO_Y(1)$ and $\SE \otimes \SO_Y(1)$ are nonspecial.
Let $\, \PP^s \subset \PP^r$ denote, respectively, the projective spaces of (one quotients) of $\,H^0(\SO_Y(1))$ and $H^0(\SO_Y(1)) \oplus H^0(\SE \otimes \SO_Y(1))$.
Assume that $r \geq 3$.
Let $Y \subset \PP^r$ be the embedding defined as the composition of the embedding $Y \subset \PP^s$ given by the complete linear series $H^0(\SO_Y(1))$ and $\PP^s \subset \PP^r$.\\
Assume that $\WY \subset \PP^r$ is a nondegenerate embedded ribbon over $Y \subset \PP^r$ with conormal bundle $\SE$.
Assume that $p_a(\WY) \geq 3$.
Assume that there is a smooth irreducible double cover $X \overset{\pi} \to Y$ with $\pi_*\SO_X/\SO_Y = \SE$.
Let $X \overset{\varphi}\to \PP^r$ be the morphism obtained as the composition of $\pi$ with the inclusion of $Y$ in $\PP^r$.
Then
\begin{enumerate}
\item there exists a smooth irreducible family $\SX$ proper and flat over a smooth pointed affine curve $(T, 0)$ and a $T$--morphism $\SX \overset{\Phi}\to \PP_T^r$ with the following properties:
\begin{enumerate}
\item the general fiber $\SX_t \overset{\Phi_t}\to \PP^r, t \neq 0,$ is a closed immersion of a smooth irreducible projective curve $\SX_t$,
\item the central fiber $\SX_0 \overset{\Phi_0}\to \PP^r$ is $X \overset{\varphi} \to \PP^r$; and
\end{enumerate}
\item the image of $\SX \overset{\Phi}\to \PP_T^r$ is a closed integral subscheme $\SY \subset \PP^r_T$ flat over $T$ with the following properties:
\begin{enumerate}
\item the general fiber $\, \SY_t, t \neq 0,$ is a smooth irreducible projective nondegenerate curve with nonspecial hyperplane section in $\PP^r$,
\item the central fiber $\, \SY_0$ is $\,\WY \subset \PP^r$.
\end{enumerate}
\end{enumerate}
\end{theorem}
We remark that the condition $p_a(\WY) \geq 3$ is imposed in Theorem~\ref{embsmoothing.0} for technical reasons regarding its proof.\\
To obtain Theorem~\ref{abssmoothing.0} from Theorem~\ref{embsmoothing.0} we need a criterion to decide whether every ribbon with a fixed conormal bundle $\SE$ over a curve $Y$ can be embedded as a nondegenerate subscheme in the same projective space $\PP^r$ extending an embedding $Y \hookrightarrow \PP^r$. This criterion is Proposition~\ref{nondegembedding}.
As a consequence we see in Theorem~\ref{embedding} how to obtain nondegenerate projective embedding in the same projective space for all ribbons with a fixed conormal bundle.\\\\
The proof of Theorem~\ref{embsmoothing.0} consists in extending an infinitesimal deformation of the map $X \overset{\varphi}\to \PP^r$ to a deformation over an affine base.
In order that the image of the deformation over the affine base contains $\WY$ as central fiber, we pick the infinitesimal deformation out so that this is already true at the infinitesimal level, i.e. we previously prove that the ribbon can be infinitesimally smoothed.
This key result on infinitesimal smoothing, which is a direct consequence of the Theorem~\ref{key-curves} obtained from the general infinitesimal theory that we develop in Section $3$, is:
\begin{theorem}\label{key.curves.0}
Let $Y$ be a smooth irreducible projective curve in $\PP^r$ and let $\SE$ be a line bundle on~$Y$.
Assume that there is a smooth irreducible double cover $X \overset{\pi} \to Y$ with $\pi_*\SO_X/\SO_Y = \SE$.
Let $X \overset{\varphi}\to \PP^r$ be the morphism obtained as the composition of $\pi$ and the inclusion of $Y$ in $\PP^r$.
Then every ribbon over $Y$, with conormal bundle $\SE$, embedded in $\PP^r$ is the central fiber of the image of some first--order infinitesimal deformation of $\varphi$.
\end{theorem}
To appreciate the scope of our results we specialize them to particular cases of $Y$, see Corollary~\ref{abssmoothing2}, Corollary~\ref{elliptic-smoothing} and Corollary~\ref{rational-smoothing}.\\
If $g=0$ some of the facts we obtain are that
\begin{enumerate}
\renewcommand{\theenumi}{\alph{enumi}}
\item all rational ribbons of arithmetic genus $h \geq 3$ are infinitesimally produced by a hyper\-elliptic curve asso\-ciated to the conormal bundle of the ribbons $\SO_{\PP^1}(- h -1)$,
\item all rational ribbons of arithmetic genus $h \geq 3$ can be embedded in $\PP^{h}$ with degree $2h$ over a rational normal curve,
\item all rational ribbons of arithmetic genus $h \geq 3$ and degree $2h$ over a rational normal curve in $\PP^{h}$ are smoothable,
\item hence all rational ribbon of arithmetic genus greater than or equal to $3$ are smoothable.
\end{enumerate}
If $g=1$, i.e. for an elliptic curve $Y$, some of the analogous facts we prove are that
\begin{enumerate}
\renewcommand{\theenumi}{\alph{enumi}}
\item all ribbons of arithmetic genus $h \geq 3$ are infinitesimally produced by bi--elliptic curves asso\-ciated to the conormal bundles of the ribbons,
\item all ribbons of arithmetic genus $h \geq 6$ can be embedded in $\PP^{h-2}$ with degree $2h-2$ over an elliptic normal curve $Y \subset \PP^{h-2}$ with hyperplane section not isomorphic to the dual of the conormal bundle of the ribbons, all ribbons of arithmetic genus $h=4$ or $5$ can be embedded in $\PP^{h}$ with degree $2h$ over an elliptic normal curve $Y \subset \PP^{h-1}$ and all ribbons of arithmetic genus $h=3$ can be embedded in $\PP^5$ with degree $8$ over an elliptic normal curve $Y \subset \PP^3$,
\item all ribbons embedded like above are smoothable,
\item hence all ribbons of arithmetic genus greater than or equal to $3$ over an elliptic curve are smoothable.\\
\end{enumerate}

The results of this article have other applications. In~\cite{GGP}, we build on the methods of this paper to prove results on smoothing of ropes of arbitrary multiplicity and apply it to study in detail ropes of multiplicity three on $\PP^1$.\\\\
{\em Conventions.}
We work over a fixed algebraically closed field $k$ of zero characteristic. All schemes considered are separated and of finite type over $k$.
\section{Preliminaries}
The definitions and facts gathered here are known in the references \cite[\S$1$]{BayerEisenbud95}, \cite[\S$1$]{GallegoPurna97}, \cite[\S$2$]{HulekVandeVen85} for ribbons.
We state them here for ropes without proofs, the ones in the references translate almost word by word.
\begin{noname}\label{noname1.1}
{\rm Let $Y$ be a reduced connected scheme and let $\SE$ be a locally free sheaf of rank $n-1$ on $Y$.
An $n$-rope over $Y$ with conormal bundle $\SE$ is a scheme $\WY$ with ${\WY}_{\mathrm{red}}=Y,$ such that $\SI_{Y, \WY}^2=0$ and $\SI_{Y, \WY} \simeq \SE$ as $\SO_{Y}$--modules.
If $\SE$ is a line bundle, $\WY$ is called a ribbon over $Y$.}
\end{noname}
\begin{noname}\label{noname1.2}
{\rm A rope $\WY$ over $Y$ with conormal bundle $\SE$ is determined by the {\em extension class} $[e_{\WY}] \in \mathrm{Ext}_Y^1(\Omega_Y, \SE)$ of its {\em restricted cotangent sequence}, the lower exact sequence in the pullback diagram
\begin{equation*}\label{cotangent}
\xymatrix@C-5pt@R-10pt{
0 \ar[r] & \SE \ar@{=}[d] \ar[r] & \SO_{\WY} \ar[d] \ar[r] & \SO_Y \ar[d] \ar[r] & 0 \\
0 \ar[r] & \SE \ar[r] & {\Omega_{\WY}|}_Y \ar[r] &  \Omega_Y \ar[r] &  0. }
\end{equation*}
If $\WY'$ is another rope over $Y$ with conormal bundle $\SE$ then $\WY$ and $\WY'$ are isomorphic over $Y$ iff its extension classes $[e_{\WY}]$ and $[e_{\WY'}]$ are in the same orbit by the action of the automorphisms of $\,\SE$ in $\mathrm{Ext}^1(\Omega_Y, \SE)$.\\
The rope associated with the split class is the unique rope $\WY$ such that the inclusion $Y \hookrightarrow \WY$ admits a retraction. This rope is called the {\em split rope}.}
\end{noname}
\begin{noname}\label{noname1.3}
{\rm An embedded rope $\WY$, with conormal bundle $\SE$, over a smooth irreducible closed subvariety $Y$ of a smooth irreducible variety $Z$ is defined by a subbundle of $\SN_{Y,Z}$ with dual bundle isomorphic to $\SE$.
The ideal of $\WY$ inside $Z$ is the kernel of the surjective composite homomorphism $ \SI_{Y,Z} \twoheadrightarrow \SI_{Y,Z}/\SI_{Y,Z}^2 \twoheadrightarrow \SE $.}
\end{noname}
\begin{noname}\label{noname1.4}
{\rm Let $Y$ be a reduced connected scheme and let $Z$ be a scheme.
The extension morphisms $\WY \overset{\wi}\to Z$, to a rope $\WY$ over $Y$ with conormal bundle $\SE$, of a given morphism $Y \overset{i}\to Z$ are in one--to--one correspondence with the homomorphisms $i^* \Omega_Z \overset{\omega}\to {\Omega_{\WY}|}_Y$ making the diagram
\begin{equation*}\label{pull-back-Di}
\xymatrix@C-5pt@R-9pt{
 &  &  & i^*\Omega_Z \ar[dl]_-{\omega} \ar[d]^-{\mathrm{D}i} &  \\
0 \ar[r] & \SE \ar[r]^-{j} & {{\Omega_{\WY}|}_Y} \ar[r]^-{p} &  \Omega_Y \ar[r] & 0,}
\end{equation*}
commutative, i.e. with the splittings of the exact sequence with class the image of $[e_{\WY}]$ by the map $\mathrm{Ext}_Y^1(\Omega_Y, \SE) \overset{\mathrm{D}i}\to \mathrm{Ext}_Y^1(i^*\Omega_Z, \SE)$.
In particular an extension morphism exists iff $\mathrm{D}i([e_{\WY}])=0$.}
\end{noname}
\begin{noname}\label{noname1.5}
{\rm Let $Y$ be a reduced connected scheme and let $Z$ be a scheme.
Let $\SE$ be a locally free sheaf of rank $n-1$ on $Y$.
The algebra of the $n$--rope $\WY$ with extension class $[e]$, where $e$ is the exact sequence $0 \to \SE \overset{j}\to \SG \overset{p}\to \Omega_Y \to 0$, is, as a sheaf of abelian groups, the pullback of the homomorphisms $\SO_Y \overset{\mathrm{d}}\to \Omega_Y$ and $\SG \overset{p}\to \Omega_Y$.
So the sections of $\SO_{\WY}$ over an open set $U \subset Y$ are 
\begin{equation*}
\SO_{\WY}=\{(c, s) \in \SO_Y \oplus \SG \, | \, \mathrm{d}c=ps\},
\end{equation*}
with $k$--algebra structure defined by $(c, s)(c', s')=(c \, c',c\,s'+c's)$.\\
Let $Y\overset{i}\to Z$ be a morphism.
The extension morphism $\WY \overset{\wi}\to Z$ associated with a commutative diagram
\begin{equation}\label{e-omega}
\xymatrix@C-5pt@R-9pt{
&&&i^*\Omega_Z \ar[dl]_-{\omega} \ar[d]^-{\mathrm{D}i}&\\
0 \ar[r] & \SE \ar[r]^-{j} & \SG \ar[r]^-{p} &  \Omega_Y \ar[r] &  0,}
\end{equation}
is defined by the unique map of $k$--algebras $\SO_Z \overset{\wi^{\sharp}}{\to} i_*\SO_{\WY}$ making the diagram
\begin{equation}\label{wi^sharp}
\xymatrix@C-5pt@R-7pt{
\SO_Z \ar[dd] \ar[dr]|-{\hole \wi^\sharp} \ar@/^/[drr]^-{i^\sharp}&&\\
&i_*\SO_{\WY} \ar[d] \ar[r]&i_*\SO_Y \ar[d]^-{i_*\mathrm{d}}\\
\Omega_Z \ar[r]^-{\omega^\flat}&i_{*}\SG \ar[r]^-{i_*p} &i_*\Omega_Y}
\end{equation}
commutative, where $\omega^\flat$ corresponds to $\omega$ by the adjunction isomorphism.\\
If $W \subset Z$ and $U \subset Y$ are open affine subsets with $U \subset i^{-1} W$ and $a$ is a section of $\SO_Z$ over $W$ then, from~\eqref{wi^sharp}, the map $\wi^\sharp$ is written in the form
\begin{equation}\label{wi^sharp-formula}
\wi^{\sharp}a = (i^{\sharp} a, \,\, \omega(\mathrm{d}a \otimes 1)).
\end{equation}}
\end{noname}
\section{Spaces parametrizing infinitesimal extensions of morphisms}
Let $Y$ be a smooth irreducible closed subvariety of a smooth irreducible variety $Z$.
Let $Y \overset{i}{\hookrightarrow} Z$ be the closed immersion.
Then we have an exact sequence on $Y$
\begin{equation*}\label{conor-Y}
\mathrm{Hom}(i^*\Omega_Z,\SE) \to \mathrm{Hom}(\SN_{Y,Z}^*,\SE) \overset{\delta}\to
\mathrm{Ext}^1(\Omega_Y, \SE) \overset{\mathrm{D}i}\to \mathrm{Ext}^1(i^*\Omega_Z, \SE).
\end{equation*}
Therefore, from~\ref{noname1.4}, a morphism $\WY \overset{\wi}\to Z\,$ extension of $Y \overset{i}\to Z\,$ exists iff $[e_{\WY}]$ admits a lifting by $\delta$ to an element $\tau \in \mathrm{Hom}(\SN_{Y,Z}^*,\SE)$.\\
The following result tells us that the extensions of $i$ to $\WY$ are in one--to--one correspondence with the liftings of $[e_{\WY}]$.
\begin{proposition}\label{extension}
Let $Y$ be a smooth irreducible closed subvariety of a smooth irreducible variety $Z$.
Let $Y \overset{i}{\hookrightarrow} Z$ be the closed immersion.
Let $\SE$ be a locally free sheaf of rank $n-1$ on $Y$.
\begin{enumerate}
\item There is a one--to--one correspondence between pairs $(\WY, \wi)$, where $\WY$ is an $n$--rope over $Y$ with conormal bundle $\,\SE$ and $\, \WY \overset{\wi}{\to} Z$ is a morphism extending $Y \overset{i}{\hookrightarrow} Z$ and elements $\tau \in \mathrm{Hom}(\SN_{Y,Z}^*,\SE)$.
Moreover, if $\tau$ and $(\WY, \wi)$ are in correspondence then $\delta \tau =[e_{\WY}]$.
Two pairs are isomorphic over $Y$ iff the corresponding elements are in the same orbit by the action of the automorphisms of $\,\SE$ in $\mathrm{Hom}(\SN_{Y,Z}^*,\SE)$.
\item The image subscheme of the map $\, \wi \,$ induced by $\tau$ has ideal in $Z$ equal to the kernel of the composite homomorphism $\SI_{Y,Z} \to \SN_{Y,Z}^* \overset{\tau} \to \SE $. Moreover, $\wi$ is a closed immersion iff $\tau$ is surjective.
\end{enumerate}
\end{proposition}
\begin{Proof}
Let $e$ be an extension $0 \to \SE \overset{j}\to \SG \overset{p}\to \Omega_Y \to 0$ and let $i^*\Omega_{Z}\overset{\omega}\to \SG$ be a homomorphism like in a commutative diagram~\eqref{e-omega}.
From \ref{noname1.4}, we obtain a one--to--one correspondence between extension pairs $(\WY, \wi)$ and classes $[(e,\omega)]$ under the obvious equivalence between pairs $(e, \omega)$.
Moreover, two extension pairs are isomorphic iff the corresponding classes are in the same orbit by the action of the automorphisms of $\SE$ in $\mathrm{Ext}^1(\Omega_Y,\SE)$.\\
To prove (1) we set up a bijection, compatible with the action of the automorphisms of $\SE$, between homomorphisms $\SN_{Y,Z}^*\overset{\tau}\to \SE$ and classes $[(e,\omega)]$, as follows:\\
Let $\SN_{Y,Z}^*\overset{\tau}\to \SE$ be a homomorphism, recall that the extension class $\delta \tau \in \mathrm{Ext}^1(\Omega_Y,\SE)$ is represented by the lower exact sequence in the push--out diagram:
\begin{equation}\label{push-tau}
\xymatrix@C-5pt@R-10pt{
0 \ar[r] & \SN_{Y,Z}^* \ar[d]_-{\tau} \ar[r]^-{j'} & i^*\Omega_Z \ar[d]^-{\omega_\tau} \ar[r]^-{\mathrm{D}i} & \Omega_Y \ar@{=}[d] \ar[r] &0\\
0 \ar[r] & \SE \ar[r]^-{j} &\frac{\SE \oplus i^*\Omega_Z}{\mathrm{im}(-\tau \oplus j')} \ar[r]^-{p} & \Omega_Y \ar[r] &0.}
\end{equation}
We denote $(i^*\Omega_Z)_\tau = \frac{\SE \oplus i^*\Omega_Z}{\mathrm{im}(-\tau \oplus j')}$ and we call $e_{\tau}$ the extension $0 \to \SE \overset{j}\to (i^*\Omega_Z)_\tau \overset{p}\to \Omega_Y \to 0$.
In this way to each $\SN_{Y,Z}^*\overset{\tau}\to \SE$ we assign the pair $(e_{\tau}, \omega_{\tau})$ defined by the commutative diagram:
\begin{equation*}\label{e-omega^tau}
\xymatrix@C-5pt@R-9pt{
&&&i^*\Omega_Z \ar[dl]_-{\omega_\tau} \ar[d]^-{\mathrm{D}i}&\\
0 \ar[r] & \SE \ar[r]^-{j} &  (i^*\Omega_Z)_\tau \ar[r]^-{p} &  \Omega_Y \ar[r] &  0. }
\end{equation*}
In the opposite direction, fix a pair $(e,\omega)$ defined by a diagram like~\eqref{e-omega}, then there is a unique homomorphism $\tau^\omega$ making the diagram
\begin{equation}\label{tau^omega}
\xymatrix@C-5pt@R-10pt{
0 \ar[r] & \SN_{Y,Z}^* \ar[d]_-{\tau^{\omega}} \ar[r]^-{j'} & i^*\Omega_Z \ar[d]^-{\omega} \ar[r]^-{\mathrm{D}i} & \Omega_Y \ar@{=}[d] \ar[r] &0\\
0 \ar[r] & \SE \ar[r]^-{j} &\SG \ar[r]^-{p} & \Omega_Y \ar[r] &0}
\end{equation}
commutative.
In this way to each pair $(e, \omega)$ we assign a homomorphism $\SN_{Y,Z}^* \overset{\tau^\omega}\longrightarrow \SE$.\\
It is easy to verify that
this establishes the bijection between homomorphisms $\SN_{Y,Z}^*\overset{\tau}\to \SE$ and classes $[(e,\omega)]$ and that this bijection is compatible with the action of the automorphisms of $\SE$ in both sets.
Therefore (1) is proved.\\
Now we prove (2).
Fix $\tau \in \mathrm{Hom}(\SN_{Y,Z}^*,\SE)$ and let $(\WY, \wi)$ be the extension pair defined by $(e_{\tau}, \omega_{\tau})$ as in~(1).
Let $\SJ$ denote the kernel of $\wi^{\sharp}$.
If we denote the composition $\SI_{Y,Z} \twoheadrightarrow \SN_{Y,Z}^* \overset{\tau}\to \SE$ also by $\SI_{Y,Z} \overset{\tau}\to \SE$ then, from the definition of $(i^*\Omega_Z)_\tau$ and~\eqref{wi^sharp-formula}, we see that the diagram
\begin{equation*}
\xymatrix@C-5pt@R-10pt{
\SI_{Y,Z} \ar[d]_-{\tau} \ar[r] & \SO_Z \ar[d]^-{\wi^{\sharp}}\\
\SE \ar[r] & \SO_{\WY}}
\end{equation*}
is commutative.
Therefore we obtain a commutative exact diagram
\begin{equation*}\label{jota-ideal}
\xymatrix@C-5pt@R-14pt{
 & 0 \ar[d] & 0 \ar[d] & &\\
 & \SJ \ar[d] \ar@{=}[r] & \SJ \ar[d] & & \\
0 \ar[r] & \SI_{Y,Z} \ar[d]_-{\tau} \ar[r] & \SO_Z \ar[d]^-{\wi^{\sharp}} \ar[r] & \SO_Y \ar@{=}[d] \ar[r] & 0 \\
0 \ar[r] & \SE \ar[r] & \SO_{\WY} \ar[r] & \SO_Y \ar[r] & 0.}
\end{equation*}
From the Snake~Lemma it follows that $\wi^{\sharp}$ is surjective iff $\tau$ is surjective.
\end{Proof}
\begin{noname}
{\rm First--order, locally trivial, infinitesimal deformations of a morphism.\\
Let $X \overset{\varphi}{\to} Z$ be a morphism, where $X$ is a reduced connected scheme and $Z$ is a scheme.
Let $\Delta$ denote $\mathrm{Spec}\,k[\epsilon]/\epsilon^2$.
We are interested in the first--order, locally trivial, infinitesimal deformations of the pair $(X,\varphi)$, i.e. $\Delta$--morphisms $\WX \overset{\wphi}{\to} Z \times \Delta$ with central fiber $X \overset{\varphi}{\to} Z$, where $\WX$ is a first--order, {\em locally trivial}, infinitesimal deformation of $X$.\\
An space classifying them is described in \cite{Horikawa74}.
Let $\mathcal{V}=(V)$ be an open affine cover of $X$.
As usual, we let $\mathcal{C}^0(\mathcal{V},-)$ and $\mathcal{Z}^1(\mathcal{V},-)$ denote, respectively, the group of 0--cochains and 1--cocycles with respect to the covering $\mathcal{V}$ and $\delta$ the coboundary map.
Let $\varphi^*\Omega_Z \overset{\mathrm{D}\varphi}\longrightarrow \Omega_X$ be the homomorphism induced by $X \overset{\varphi}{\to} Z$.
Set
\begin{equation}\label{definitionDef}
D(X,\varphi)= \frac{\{(g,\rho) \in \mathcal{C}^0(\mathcal{V},\mathcal{H}om_{\SO_X}(\varphi^* \Omega_Z,\SO_X)) \times \mathcal{Z}^1(\mathcal{V},\mathcal{H}om_{\SO_X}(\Omega_X,\SO_X)) \; | \; \delta g = \rho \, \mathrm{D}\varphi \}}{\{(h \, \mathrm{D}\varphi,\delta h) \; | \; h \in \mathcal{C}^0(\mathcal{V},\mathcal{H}om_{\SO_X}(\Omega_X,\SO_X))\}}.
\end{equation}}
\end{noname}
\begin{lemma}\label{Def}{\rm \cite[4.2]{Horikawa74}}
Let $X \overset{\varphi}{\to} Z$ be a morphism, where $X$ is a reduced connected scheme and $Z$ is a scheme.
Let $D(X,\varphi)$ defined by~\eqref{definitionDef}.
Then\\
{\rm (1)} $D(X,\varphi) \,$ does not depend on the affine cover.\\
{\rm (2)} We have two exact sequences
\begin{eqnarray*}\label{exactDef}
\begin{aligned}
&{\rm Hom}(\Omega_X,\SO_X) \overset{\mathrm{d}\varphi} \to 
{\rm Hom}(\varphi^*\Omega_Z,\SO_X) \to D(X,\varphi) \to H^1(\mathcal{H}om(\Omega_X,\SO_X)) \overset{\mathrm{d}\varphi} \to 
H^1(\mathcal{H}om(\varphi^*\Omega_Z,\SO_X)),\\
&0 \to  H^1(\mathcal{H}om(\Omega_{X/Z}, \SO_X)) \to D(X,\varphi) \to H^0(\mathcal{N}_{\varphi}) \to H^2(\mathcal{H}om(\Omega_{X/Z}, \SO_X)),
\end{aligned}
\end{eqnarray*}
where $\mathcal{N}_{\varphi}$ denotes the cokernel of $\,\mathcal{H}om(\Omega_X, \SO_X) \overset{\mathrm{d}\varphi}{\longrightarrow}\mathcal{H}om(\varphi^*\Omega_Z, \SO_X)$.\\
In particular, if the map 
$\mathcal{H}om(\Omega_X, \SO_X)\overset{\mathrm{d}\varphi}{\longrightarrow}
\mathcal{H}om(\varphi^*\Omega_Z, \SO_X)$ is injective then there is a natural isomorphism $D(X, \varphi) \simeq H^0(\mathcal{N}_{\varphi})$.
\end{lemma}
\begin{remark}\label{Def2}
{\rm The content of Lemma~\ref{Def} is a purely cohomological fact of sheaves on $X$. The proof given in \cite[4.2]{Horikawa74}, in the context of complex manifolds, is still valid in the case of a homomorphism of quasi--coherent sheaves $\SA \overset{F}\to\SB$ with kernel and cokernel, respectively, $\SK$ and $\SN$, on a noetherian separated scheme $X$. Let $\SV$ be an open cover of $X$.
We set
\begin{equation}\label{definitionDef2-V}
D(X, F; \SV)= \frac{\{(g,\rho) \in \mathcal{C}^0(\mathcal{V},\SB) \times \mathcal{Z}^1(\mathcal{V},\SA) \; | \; \delta g = F \rho \}}{\{(F h,\delta h) \; | \; h \in \mathcal{C}^0(\mathcal{V},\SA)\}}
\end{equation}
and 
\begin{equation}\label{definitionDef2}
D(X, F)= \underset{\underset{\SV}\longrightarrow}\lim \,\,D(X,F;\SV),
\end{equation}
where the direct limit is taken under refinement of coverings.
Then for every open {\em affine} cover $\SV$ the natural map $D(X,F;\SV) \to D(X,F)$ is an isomorphism.
Two exacts sequences, like in Lemma~\ref{Def}, are then obtained
\begin{equation*}
\begin{aligned}
& H^0(\SA) \to H^0(\SB) \to D(X, F)\to H^1(\SA) \to H^1(\SB),\\
0 \to & H^1(\SK) \to D(X, F) \to H^0(\SN)\to H^2(\SK).
\end{aligned}
\end{equation*}
{\hfill $\square$}
}
\end{remark}
\begin{proposition}\label{extension2}
Let $X \overset{\varphi}{\to} Z$ be a morphism, where $X$ is a reduced connected scheme and $Z$ is a scheme.
Let $D(X,\varphi)$ be defined by~\eqref{definitionDef}. There is a one--to--one correspondence between pairs $(\WX, \widetilde{\varphi})$ up to $\Delta$--isomorphism, where $\WX$ is a first--order, locally trivial, infinitesimal deformation of $\,X$ and $\WX \overset{\widetilde{\varphi}}\to Z \times \Delta$ is a $\Delta$--morphism with central fiber $X \overset{\varphi}\to Z$, and classes in $D(X, \varphi)$.
\end{proposition}
\begin{Proof}
A first--order infinitesimal deformation $\WX$ is a ribbon over $X$ with conormal $\SO_X$ and conversely.
The $\Delta$--morphisms $\WX \overset{\wphi}{\to} Z \times \Delta$ are in bijection with the $k$--morphisms $\WX \overset{\whphi}{\to} Z$ and 
a $\Delta$--morphism $\wphi$ is an extension of $\varphi$ iff the corresponding $k$--morphism $\whphi$ is an extension of $\varphi$.\\
Now, arguing like in the proof of Proposition~\ref{extension}, we establish a bijection between extension pairs $(\WX, \wphi)$ up to $\Delta$--isomorphism and classes of pairs $[(e,\omega)]$ defined by equivalence of diagrams
\begin{equation}\label{e-omega-phi}
\xymatrix@C-5pt@R-10pt{
&&&\varphi^*\Omega_Z \ar[dl]_-{\omega} \ar[d]^-{\mathrm{D}\varphi}&\\
0 \ar[r] & \SO_X \ar[r]^-{j} & \SG \ar[r]^<(.4){p} &  \Omega_X \ar[r] &  0, }
\end{equation}
by assigning to $(\WX, \wphi)$ the pair $(e,\omega)$ associated with $(\WX, \whphi)$.
We observe, for future reference, that the correspondence $\wphi \leftrightarrow \whphi$ is locally expressed by
\begin{equation}\label{local-wphi-whphi}
\wphi^\sharp(a+a'\epsilon)=\whphi^\sharp(a)+\whphi^\sharp(a')\epsilon, \; \mathrm{where}\; a+a'\epsilon \in \SO_Z \oplus \SO_Z\epsilon.
\end{equation}
From~\eqref{local-wphi-whphi} and~\eqref{wi^sharp-formula} we see that the $\Delta$--morphism $\wphi$ associated with $(e,\omega)$ is locally expressed by:
\begin{equation}\label{phi-tilde-local}
\begin{aligned}
\SO_Z \oplus \SO_Z\epsilon & \overset{\wphi^{\sharp}}{\longrightarrow} \varphi_*\SO_{\WX}\\
a+a'\epsilon & \mapsto (\varphi^{\sharp}a,\, \omega(\mathrm{d}a \otimes 1)+j\varphi^{\sharp}a'),
\end{aligned}
\end{equation}
where
$\SO_{\WX}=\{(b,s) \in \SO_X \oplus \SG \, | \, p\,s=\mathrm{d}b\}$.
Now, to prove Proposition~\ref{extension2}, we establish a bijection between the set of classes $[(e,\omega)]$, where the extension $e$ is locally split, and $D(X,\varphi)$.\\
Locally split extension classes are classified by the subspace $H^1(\mathcal{H}om_{\SO_X}(\Omega_X,\SO_X))  \subset \mathrm{Ext}^1(\Omega_X,\SO_X)$.
Recall that the inclusion map, obtained from the spectral sequence of local and global Ext's, is as follows: to a class $[\rho]$, with $\rho \in \mathcal{Z}^1(\mathcal{V},\mathcal{H}om_{\SO_X}(\Omega_X,\SO_X))$, corresponds the class of the extension $e$ defined by the gluing diagrams
\begin{equation}\label{glue-G}
\xymatrix@C-5pt@R-10pt{
0 \ar[r] & \SO_{\smV \cap \smV'} \ar@{=}[d] \ar[r] & \SO_{\smV\cap \smV'} \oplus {\Omega_{\smV\cap \smV'}} \ar[d]^-{\sigma_{\smV \smV'}} \ar[r] & {\Omega_{\smV\cap \smV'}} \ar@{=}[d] \ar[r] & 0 \\
0 \ar[r] & \SO_{\smV\cap \smV'} \ar[r] & \SO_{\smV\cap \smV'} \oplus {\Omega_{\smV\cap \smV'}} \ar[r] & {\Omega_{\smV\cap \smV'}}  \ar[r] & 0,}
\end{equation}
where
\begin{equation*}
\sigma_{\smV \smV'}=\begin{bmatrix}\mathrm{id}_{\SO} & {\rho}_{\smV \smV'}\\ 0 & \mathrm{id}_{\Omega}\end{bmatrix}.
\end{equation*}
We set up the bijection.
Start with a class $[(e,\omega)]$, where $e$ is locally split.
Observe that a locally split extension of coherent sheaves on a noetherian separated scheme is split over every open affine subset.
Take an open affine cover $\mathcal{V}=(V)$ of $X$ so that $e$ is split in every open set of $\SV$, and take a family of local retractions $r=(r_{\smV})_{\smV\in \mathcal{V}}$ for $j$.
Then there is a unique $\rho \in \mathcal{Z}^1(\mathcal{V},\mathcal{H}om_{\SO_X}(\Omega_X,\SO_X))$ such that $\delta  r = \rho  p$.
Moreover, the cochain $r\,\omega$ verifies that $\delta(r \,\omega)=\rho\, \mathrm{D}\varphi$. Therefore the pair $(r\omega, \rho)$ defines a class in $D(X,\varphi)$.
Conversely, to a class $[(g,\rho)]$ we assign the class of the pair $(e,\omega)$ constructed as follows:
the extension $e$ is defined from $\rho$ by means of the gluing diagrams~\eqref{glue-G}.
The condition $\delta g= \rho \, \mathrm{D}\varphi $, in the definition of $D(X,\varphi)$, implies that the homomorphisms $ \varphi^*{\Omega_Z |}_{\smV} \overset{({g}_{{}_{\smV}},\mathrm{D}\varphi)}\longrightarrow \SO_{\smV} \oplus {\Omega_{\smV}}$, glue to define a global homomorphism $\omega$ and a commutative diagram like~\eqref{e-omega-phi}.
It is easy to verify that this establishes the bijection.
\end{Proof}
\section{Images of first-order infinitesimal deformations of certain type of morphisms}
In this section we set up the general process relating ropes, over a smooth irreducible closed subvariety $Y$ of a smooth irreducible variety $Z$, mapping to the ambient variety $Z$ and first--order, locally trivial, infinitesimal deformations of morphisms which are the composition of a finite cover of $Y$ and the inclusion $Y \hookrightarrow Z$.\\
We start by fixing the setting for this section:
\begin{noname}\label{setting}
{\rm Let $X\overset{\varphi}\to Z$ be a morphism from an integral Cohen--Macaulay variety $X$ to a smooth irreducible variety $Z$.
Let $Y$ be the (scheme--theoretic) image of $\varphi$.
Let $Y \overset{i}\hookrightarrow Z$ denote the closed immersion.
Assume that $Y$ is smooth and that $\varphi$ induces a finite morphism $X \overset{\pi}\to Y$.\\
In this conditions $\pi$ is surjective and flat.
The algebra $\pi_*\SO_X$ is a locally free $\SO_Y$--module of some rank $n$ and the trace map gives a splitting for the injective map $\SO_Y \to \pi_{*}\SO_X$.
Hence $\pi_{*}\SO_X$ is the direct sum of $\SO_Y$ and a rank $n-1$ locally free $\SO_Y$--module $\SE$.}
\end{noname}
\begin{noname}\label{normal-varphi}
{\rm Let $\SN_{\varphi}$ be the normal sheaf of the morphism $X \overset{\varphi}\to Z$ defined by the exact sequence
\begin{equation*}\label{normal-phi}
\mathcal{H}om_{\SO_X}(\Omega_X,\SO_X) \overset{\mathrm{d}\varphi}{\longrightarrow}
\mathcal{H}om_{\SO_X}(\varphi^* \Omega_Z,\SO_X)\to \SN_{\varphi} \to 0.
\end{equation*}}
\end{noname}
This sheaf fits into a useful extension:
\begin{lemma}\label{exact-normals}
In the conditions of~\ref{setting} there is an exact sequence
\begin{equation}\label{normal-extension}
0 \to \SN_{\pi} \to \SN_{\varphi} \to \pi^* \SN_{Y,Z} \to 0.
\end{equation}
\end{lemma}
\begin{Proof}
We have the sequence exact also on the left 
\begin{equation*}\label{cotangentX-Y}
\xymatrix@C-5pt{
0 \ar[r]& \pi^*\Omega_Y \ar[r]^-{\mathrm{D}\pi}& \Omega_X \ar[r]& \Omega_{X/Y} \ar[r]& 0.}
\end{equation*}
So we obtain the exact sequence
\begin{equation*}\label{dif-pi-dual}
\xymatrix@C-5pt{
0 \ar[r] & \mathcal{H}om_{\SO_X}(\Omega_X, \SO_X) 
\ar[r]^-{\mathrm{d}\pi} & \mathcal{H}om_{\SO_X}(\pi^*\Omega_Y, \SO_X) 
\ar[r] & \SN_{\pi} \ar[r] & 0.}
\end{equation*}
From $\,0 \to \mathcal{I}_{Y,Z}/\mathcal{I}_{Y,Z}^2 \to i^*\Omega_Z \overset{\mathrm{D}i}\longrightarrow \Omega_Y \to 0\,$ we obtain the exact sequence
\begin{equation*}
0 \to \mathcal{H}om_{\SO_X}(\pi^*\Omega_Y, \SO_X) \overset{\pi^*\mathrm{d}i}\longrightarrow \mathcal{H}om_{\SO_X}(\varphi^*\Omega_Z, \SO_X) \to \mathcal{H}om_{\SO_X}(\pi^*\mathcal{I}_{Y,Z}/\mathcal{I}_{Y,Z}^2,
\SO_X)\to 0.
\end{equation*}
So we see that also the map $\mathcal{H}om_{\SO_X}(\Omega_X, \SO_X) \overset{\mathrm{d}\varphi} \longrightarrow \mathcal{H}om_{\SO_X}(\varphi^*\Omega_Z, \SO_X)$ is injective.
Therefore we obtain an exact commutative diagram
\begin{equation}\label{normal-diagram}
\xymatrix@C-12pt@R-10pt{
& &   0 \ar[d] & 0 \ar[d] &  \\
0 \ar[r] & \mathcal{H}om_{\SO_X}(\Omega_X, \SO_X) \ar@{=}[d]
\ar[r]^-{\mathrm{d}\pi} & \mathcal{H}om_{\SO_X}(\pi^*\Omega_Y, \SO_X) \ar[d]^-{\pi^*\mathrm{d}i}
\ar[r] & \SN_{\pi} \ar[d] \ar[r] & 0 \\
0 \ar[r] & \mathcal{H}om_{\SO_X}(\Omega_X, \SO_X) \ar[r]^-{\mathrm{d}\varphi} &
\mathcal{H}om_{\SO_X}(\varphi^*\Omega_Z, \SO_X) \ar[d] \ar[r] &
\mathcal{N}_{\varphi} \ar[d] \ar[r]& 0 \\
 &  & \mathcal{H}om_{\SO_X}(\pi^*\mathcal{I}_{Y,Z}/\mathcal{I}_{Y,Z}^2,
\SO_X) \ar[d] \ar@{=}[r] & \pi^*\SN_{Y,Z} \ar[d]& \\
 & &  0 & 0 &}
\end{equation}
and the right--hand side column is the desired exact sequence.
\end{Proof}\\
With the notations of \eqref{definitionDef}, \eqref{definitionDef2} and Lemma~\ref{exact-normals} we have
\begin{lemma}\label{key-map} In the conditions of~\ref{setting} there is a commutative diagram
\begin{equation}\label{g-restricted}
\xymatrix@C-5pt@R-9pt{
D(X, \varphi) \ar[d]_-{\mu} \ar[r]^-{\sim} & H^0(\SN_{\varphi}) \ar[d]^-{\phi}\\
D(X, \pi^*\mathrm{d}i) \ar[d]_-{\alpha'}^-{\wr} \ar[r]^-{\sim} & \mathrm{Hom}(\pi^*\mathcal{I}/\mathcal{I}^2,\SO_X)\ar[d]^-{\alpha}_-{\wr}\\
D(Y, (\mathrm{D}i)')  \ar[r]^-{\sim} & \mathrm{Hom}(\mathcal{I}/\mathcal{I}^2,\pi_*\SO_X),}
\end{equation}
where $H^0(\SN_{\varphi}) \overset{\phi}\to \mathrm{Hom}(\pi^*\mathcal{I}/\mathcal{I}^2,\SO_X)$ is the map in global sections obtained from \eqref{normal-extension}, the composition $D(X, \varphi) \to \mathrm{Hom}(\pi^*\mathcal{I}/\mathcal{I}^2,\SO_X)$ is $[(g,\rho)] \mapsto {g_|}_{\pi^*\mathcal{I}/\mathcal{I}^2}$ and if $\alpha'\mu([(g,\rho)])=[(f,\varrho)]$ then the composition $D(X, \varphi) \to \mathrm{Hom}(\mathcal{I}/\mathcal{I}^2,\pi_*\SO_X)$ is $[(g,\rho)] \mapsto {f_|}_{\mathcal{I}/\mathcal{I}^2}$.
\end{lemma}
\begin{Proof}
The spaces $D(X, \varphi)$ and $D(X, \pi^*\mathrm{d}i)$ are associated, from Remark~\ref{Def2}, respectively, to the homomorphisms $\mathrm{d}\varphi$ and $\pi^*\mathrm{d}i$ in \eqref{normal-diagram}.
The space $D(Y, (\mathrm{D}i)')$ is associated from Remark~\ref{Def2} to the homomorphism $(\mathrm{D}i)'$ in the exact sequence
\begin{equation}\label{Di'}
0 \to \mathcal{H}om_{\SO_Y}(\Omega_Y, \pi_*\SO_X)  \overset{(\mathrm{D}i)'}\longrightarrow
\mathcal{H}om_{\SO_Y}(i^*\Omega_Z, \pi_*\SO_X)\to \mathcal{H}om_{\SO_Y}(\mathcal{I}_{Y,Z}/\mathcal{I}_{Y,Z}^2,
\pi_*\SO_X)\to 0.
\end{equation}
So if $\SU=(U)$ is an open affine cover of $Y$ and $\SV=(V=\pi^{-1}U)$ is the induced open affine cover of $X$ then $D(X,\varphi)$ is given by~\eqref{definitionDef} and $D(X, \pi^*\mathrm{d}i), D(Y, (\mathrm{D}i)')$ are obtained from~\eqref{definitionDef2-V}.
The homomorphisms $\mathrm{d}\varphi$, $\pi^*\mathrm{d}i$ in \eqref{normal-diagram} and $(\mathrm{D}i)'$ in~\eqref{Di'} are injective.
Therefore, from Remark~\ref{Def2}, we have natural isomorphisms 
$D(X, \varphi) \overset{\sim}\to H^0(\SN_{\varphi})$, $\, D(X, \pi^*\mathrm{d}i) \overset{\sim}\to \mathrm{Hom}(\pi^*\mathcal{I}/\mathcal{I}^2,\SO_X)$ and $D(Y, (\mathrm{D}i)') \overset{\sim}\to \mathrm{Hom}(\mathcal{I}/\mathcal{I}^2,
\pi_*\SO_X)$.\\
From \eqref{normal-diagram} and the functorial character of $D(X,-)$ in the exact sequences of Remark~\ref{Def2} we see that there is a unique map $D(X, \varphi)\overset{\mu}\to D(X, \pi^*\mathrm{d}i)$ making the upper square in diagram~\eqref{g-restricted} commutative.
Hence this map is $[(g, \rho)] \overset{\mu}\mapsto [(g, \rho \,\mathrm{D}\pi)]$.\\
From the adjunction isomorphism 
we obtain an isomorphism
\begin{equation}\label{Des}
D(X, \pi^*\mathrm{d}i)  \overset{\alpha'}{\underset{\sim}\longrightarrow} D(Y, (\mathrm{D}i)') \;\; \mathrm{given} \, \mathrm{by} \;\;
[(g,\rho')]  \mapsto [(f,\varrho)].
\end{equation}
Moreover, the isomorphism $D(X, \pi^*\mathrm{d}i) \overset{\sim}\to \mathrm{Hom}(\pi^*\mathcal{I}/\mathcal{I}^2,\SO_X)$ given by Remark~\ref{Def2} is defined by $[(g,\rho')] \mapsto {g_|}_{\pi^*\mathcal{I}/\mathcal{I}^2}$, i.e., for the restriction of the cochain $g$ to $\pi^*\mathcal{I}/\mathcal{I}^2$ we have $\delta ({g_|}_{\pi^*\mathcal{I}/\mathcal{I}^2}) =0$ and therefore there is a global section ${g_|}_{\pi^*\mathcal{I}/\mathcal{I}^2} \in \mathrm{Hom}(\pi^*\mathcal{I}/\mathcal{I}^2, \SO_X)$.
Alike, the isomorphism $D(Y, (\mathrm{D}i)') \overset{\sim}\to \mathrm{Hom}(\mathcal{I}/\mathcal{I}^2, \pi_*\SO_X)$ is given by $[(f,\varrho)] \mapsto {f_|}_{\mathcal{I}/\mathcal{I}^2}$.\\
Now, since $g$ corresponds to $f$ by the adjunction isomorphism, we see that ${g_|}_{\pi^*\mathcal{I}/\mathcal{I}^2}$ corresponds to ${f_|}_{\mathcal{I}/\mathcal{I}^2}$ by the isomorphism $\alpha$ in~\eqref{g-restricted}.
So we have commutativity in the lower square of~\eqref{g-restricted} and the final assertions in the Lemma follow as well.
\end{Proof}
\begin{proposition}\label{squares}
In the conditions of~\ref{setting} there is a commutative diagram
\begin{equation}\label{deltas-square}
\xymatrix@C-5pt@R-7pt{
H^0(\SN_{\varphi}) \ar[d]_-{\phi} \ar[r]^-{\delta_1} & \mathrm{Ext}^1(\Omega_X,\SO_X) \ar[d]^-{\mathrm{d}\pi}\\
\mathrm{Hom}(\pi^*\mathcal{I}/\mathcal{I}^2,\SO_X)\ar_-{\alpha}^-{\wr}[d]\ar[r]^-{\delta_2}&
\mathrm{Ext}^1(\pi^*\Omega_Y, \SO_X)\ar^-{\beta}_-{\wr}[d]\\
\mathrm{Hom}(\mathcal{I}/\mathcal{I}^2,\pi_*\SO_X)\ar[r]^-{\delta_3} &\mathrm{Ext}^1(\Omega_Y,\pi_*\SO_X).}
\end{equation}
The map $\delta_1$ sends the section $\nu \in H^0(\SN_{\varphi})$ which, from the isomorphism $D(X,\varphi) \overset{\sim}\to H^0(\SN_{\varphi})$ and Proposition~\ref{extension2}, corresponds to a first--order, locally trivial, infinitesimal deformation $(\WX,\wphi)$ of $\varphi$, to the class of the extension defining $\WX$.
The maps $\delta_3$ and $\delta_2$ are, respectively, the connecting homomorphisms obtained from the sequence $0 \to \SI/\SI^2 \to i^*\Omega_Z \to \Omega_Y \to 0$ and its pullback to $X$.
\end{proposition}
\begin{Proof}
We define $\delta_1$ as the composition of the connecting map $H^0(\SN_{\varphi}) \to H^1(\mathcal{H}om(\Omega_X,\SO_X))$, obtained from the middle sequence in~\eqref{normal-diagram}, and the natural inclusion $H^1(\mathcal{H}om(\Omega_X,\SO_X)) \hookrightarrow \mathrm{Ext}^1(\Omega_X,\SO_X)$.
If $\nu \in H^0(\SN_{\varphi})$ corresponds to $[(g, \rho)] \in D(X,\varphi)$ by the isomorphism $D(X,\varphi) \overset{\sim}\to H^0(\SN_{\varphi})$ then $H^0(\SN_{\varphi}) \to H^1(\mathcal{H}om(\Omega_X,\SO_X))$ sends $\nu$ to $[\rho]$.
So $\delta_1$ sends $\nu$ to the class of the extension defined by the gluing diagrams~\eqref{glue-G}.
Therefore, from Proposition~\ref{extension2}, we see that $\delta_1$ sends $\nu$ to the class of the extension defining the first--order, locally trivial, infinitesimal deformation $\WX$ in the pair $(\WX, \wphi)$ associated with $[(g, \rho)]$.\\
Now we prove the commutativity of the upper square in~\eqref{deltas-square}.\\
The map $\delta_2$ is the composition of the connecting map $H^0(\mathcal{H}om(\pi^*\SI/\SI^2, \SO_X)) \to$ $H^1(\mathcal{H}om(\pi^*\Omega_Y, \SO_X))$ obtained from the middle column of~\eqref{normal-diagram}, and the natural inclusion $H^1(\mathcal{H}om(\pi^*\Omega_Y, \SO_X)) \hookrightarrow \mathrm{Ext}^1(\pi^*\Omega_Y,\SO_X)$.
Furthermore, this inclusion is an isomorphism for $Y$ is smooth.
Therefore, from the upper commutative square in~\eqref{g-restricted}, the commutativity of the upper square in~\eqref{deltas-square} is equivalent to the commutativity of the diagram
\begin{equation}\label{obvious1}
\xymatrix@C-5pt@R-7pt{
D(X, \varphi) \ar[d]_-{\mu} \ar[r] & H^1(\mathcal{H}om(\Omega_X,\SO_X)) \ar[d]^-{\mathrm{d}\pi}\\
D(X, \pi^*\mathrm{d}i) \ar[r] & H^1(\mathcal{H}om(\pi^*\Omega_Y,\SO_X)),}
\end{equation}
where the horizontal maps are given by Remark~\ref{Def2} and $\mu$ is the map of~\eqref{g-restricted}.
The commutativity of~\eqref{obvious1} follows directly from the definitions.\\
Now we prove commutativity in the lower square of~\eqref{deltas-square}.
The vertical map $\alpha$ is the adjunction isomorphism.
The vertical isomorphism $\beta$ follows from the facts that $\pi$ is an affine morphism and $Y$ is smooth.
Indeed, $\mathrm{Ext}^1(\pi^*\Omega_Y,\SO_X)\simeq
H^1(\mathcal{H}om(\pi^*\Omega_Y,\SO_X))$, for $\pi^*\Omega_Y$ is locally free.
From the fact that $\pi$ is affine $H^1(\mathcal{H}om(\pi^*\Omega_Y,\SO_X))\simeq H^1(\pi_*\mathcal{H}om(\pi^*\Omega_Y,\SO_X))$.
Moreover, from the adjunction isomorphism $\pi_*\mathcal{H}om(\pi^*\Omega_Y,\SO_X)\simeq
\mathcal{H}om(\Omega_Y,\pi_*\SO_X)$.
Therefore there is an isomorphism
\begin{equation*}\label{Hom-beta'}
\xymatrix@C-7pt{
H^1(\mathcal{H}om(\pi^*\Omega_Y,\SO_X)) \ar[r]^-{\beta'}_-{\sim}&
H^1(\mathcal{H}om(\Omega_Y,\pi_*\SO_X)).}
\end{equation*}
Finally $H^1(\mathcal{H}om(\Omega_Y,\pi_*\SO_X)) \simeq
\mathrm{Ext}^1(\Omega_Y,\pi_*\SO_X),$ for $\Omega_Y$ is locally free.\\
Furthermore, the commutativity of the lower square of~\eqref{deltas-square} is equivalent to the commutativity of the diagram
\begin{equation}\label{deltas'-square}
\xymatrix@C-5pt@R-7pt{
\mathrm{Hom}(\pi^*\mathcal{I}/\mathcal{I}^2,\SO_X)\ar_-{\alpha}^-{\wr}[d]\ar[r]^-{\delta_2'}&
H^1(\mathcal{H}om(\pi^*\Omega_Y, \SO_X))\ar^-{\beta'}_-{\wr}[d]\\
\mathrm{Hom}(\mathcal{I}/\mathcal{I}^2,\pi_*\SO_X)\ar[r]^-{\delta_3'} &H^1(\mathcal{H}om(\Omega_Y,\pi_*\SO_X)),}
\end{equation}
where $\delta_2'$ and $\delta_3'$ are obtained, respectively, from the middle column of~\eqref{normal-diagram} and the exact sequence~\eqref{Di'}.
So we are reduced to prove commutativity in~\eqref{deltas'-square}.
Indeed, from~\eqref{Des} we obtain a commutative diagram
\begin{equation}\label{DX-DY}
\xymatrix@C-5pt@R-7pt{
D(X, \pi^*\mathrm{d}i)\ar[d]_-{\alpha'}^-{\wr} \ar[r]&
H^1(\mathcal{H}om(\pi^*\Omega_Y, \SO_X))\ar[d]^-{\beta'}_-{\wr}\\
D(Y, (\mathrm{D}i)')\ar[r]& H^1(\mathcal{H}om(\Omega_Y,\pi_*\SO_X)),}
\end{equation}
where the horizontal maps are given by Remark~\ref{Def2}.
Now, from the lower commutative square in~\eqref{g-restricted}, the commutativity of~\eqref{DX-DY} implies commutativity in~\eqref{deltas'-square}.
\end{Proof}\\
The next Proposition~\ref{construction} provides geometric meaning to Proposition~\ref{squares}.
It is the starting point for Theorem~\ref{image} and Theorem~\ref{main}, the main results in this section.
In order to state Proposition~\ref{construction} we need the following:
\begin{lemma}\label{glue-ropes}
Let $Y$ be a reduced connected scheme.
\begin{enumerate}
\item A rope $\WY$ over $Y$ with conormal bundle $\SE$ has a nontrivial structure of $\,\Delta$--scheme extending the structure of $\,\WY$ as $k$--scheme iff $\SE$ has a nonzero global section.
\item Let $\WY_1$ be an $n_1$--rope and let $\WY_2$ be an $n_2$--rope over $Y$ with conormal bundles, respectively, the locally free sheaves $\SE_1$ and $\SE_2$ on $Y$.
Let $\WY_1 \underset{\scriptscriptstyle{Y}}\cup \WY_2$ denote the scheme obtained by gluing of $\,\WY_1$ and $\WY_2$ along $Y$.
    \begin{enumerate}
    \item  The scheme $\WY_1 \underset{\scriptscriptstyle{Y}}\cup \WY_2$ is an $(n_1+n_2-1)$--rope over $Y$ with conormal bundle $\SE_1 \oplus \SE_2$.
    \item The natural isomorphism $\mathrm{Ext}^1(\Omega_Y,\SE_1\oplus\SE_2) \overset{\sim}\to \mathrm{Ext}^1(\Omega_Y,\SE_1)\oplus \mathrm{Ext}^1(\Omega_Y,\SE_2)$ sends the extension class of $\,\WY_1 \underset{\scriptscriptstyle{Y}}\cup \WY_2$ to the extension classes of $\,\WY_1$ and $\WY_2$.
    \item Let given a scheme $Z$ and a morphism $Y \to Z$.
    To give morphisms $\WY_1 \to Z$ and $\WY_2 \to Z$ extending the given morphism $Y \to Z$ is equivalent to give a morphism $\WY_1 \underset{\scriptscriptstyle{Y}}\cup \WY_2 \to Z$ extending $Y \to Z$.
    \end{enumerate}
\end{enumerate}
\end{lemma}
\begin{Proof} 
The $k$--scheme $\WY_1 \underset{\scriptscriptstyle{Y}}\cup \WY_2$ is $(Y, \SO)$, where $\SO$ is the subsheaf of $k$--algebras of $\SO_{\WY_1} \oplus \SO_{\WY_2}$ defined, as a sheaf of abelian groups, as the pullback of $\SO_{\WY_1} \twoheadrightarrow \SO_Y$ and $\SO_{\WY_2} \twoheadrightarrow \SO_Y$.
So there is an exact sequence
\begin{equation*}
\xymatrix@C-5pt{
0 \ar[r] & \SE_1\oplus\SE_2 \ar[r] &\SO \ar[r] & \SO_Y \ar[r] & 0,}
\end{equation*}
showing that $\WY_1 \underset{\scriptscriptstyle{Y}}\cup \WY_2$ is a rope over $Y$ with conormal bundle $\SE_1\oplus\SE_2$.
\end{Proof}
\begin{proposition}\label{construction}
Assume the conditions of~\ref{setting}.
\begin{enumerate}
\item There is a commutative diagram
\begin{equation}\label{construction1}
\xymatrix@C-5pt@R-7pt{
H^0(\SN_{\varphi}) \ar[d]_-{\Phi_1 \oplus \Phi_2} \ar[r]^-{\delta_1} & \mathrm{Ext}^1(\Omega_X,\SO_X) \ar[d]^-{\Psi_1 \oplus \Psi_2}\\
\mathrm{Hom}(\mathcal{I}/\mathcal{I}^2,\SO_Y)\oplus \mathrm{Hom}(\mathcal{I}/\mathcal{I}^2,\SE)\ar[r]^-{\delta} & 
\mathrm{Ext}^1(\Omega_Y,\SO_Y) \oplus \mathrm{Ext}^1(\Omega_Y,\SE).}
\end{equation}
Associated with every first--order, locally trivial, infinitesimal deformation $(\WX, \wphi)$ of $\, X\overset{\varphi}\to Z$, there are a pair $(\WY, \wi)$, where $\WY$ is an $n$--rope over $Y$ with conormal bundle $\SE$ and $\, \WY \overset{\wi}\to Z$ is a morphism extending $Y \overset{i}\hookrightarrow Z$, and a first--order deformation $\bar{Y} \subset Z \times \Delta$ of $\,Y$ in $Z$, in the following way: if $(\WX, \wphi)$ corresponds to $\nu \in H^0(\SN_{\varphi})$ then $\bar{Y}$ and $(\WY, \wi)$ are defined, respectively, by $\Phi_1 \nu$ and $\Phi_2 \nu$.
Furthermore, the extension classes of $\, \WX$, $\bar{Y}$ and $\WY$ are given, respectively, by $\delta_1 \nu$, $\Psi_1 \delta_1 \nu$ and $\Psi_2 \delta_1 \nu$.
\item Let $\bar{Y}\underset{\scriptscriptstyle{Y}}{\cup}\WY$ be the $(n+1)$--rope over $Y$ with conormal bundle $\SO_{Y} \oplus \SE$ obtained by gluing of $\,\bar{Y}$ and $\WY$ along $Y$.
There is a $\Delta$--morphism $\WX \overset{\widetilde{\pi}}\to \bar{Y}\underset{\scriptscriptstyle{Y}}{\cup}\WY$ extending $X \overset{\pi}\to Y$.
If $\widetilde{\psi}$ is the composition $\WX \overset{\widetilde{\pi}}\to \bar{Y}\underset{\scriptscriptstyle{Y}}{\cup}\WY \overset{\widetilde{\iota}}\to Z \times \Delta$, where $\bar{Y}\underset{\scriptscriptstyle{Y}}{\cup}\WY  \overset{\widetilde{\iota}}\to Z \times \Delta$ is
the unique $\Delta$--morphism extending both $\WY \overset{\wi} \to Z$ and $\bar{Y} \hookrightarrow Z \times \Delta$, then $\widetilde{\psi}$ is a first--order, locally trivial, infinitesimal deformation of $\varphi$ such that $\widetilde{\psi}$ and $\wphi$ are equal over the open subscheme of $\WX$ supported on the complementary of the support of $\, \Omega_{X/Y}$ and such that the difference of $\widetilde{\psi}$ and $\wphi$ is a first--order, locally trivial, infinitesimal deformation with trivial source of $X \overset{\pi}\to Y$.  
\item The image subscheme of $\, \bar{Y}\underset{\scriptscriptstyle{Y}}{\cup}\WY  \overset{\widetilde{\iota}}\to Z \times \Delta$ is the scheme-theoretic union of $\, \bar{Y}$ and the image subscheme of $\, \WY \overset{\wi} \to Z$.
\end{enumerate}
\end{proposition}
\begin{Proof}
From diagram~\eqref{deltas-square} and the splitting $\pi_*\SO_X = \SO_Y \oplus \SE$, we obtain the commutative diagram~\eqref{construction1}.
From Proposition~\ref{extension2} and Lemma~\ref{key-map} to $(\WX, \wphi)$ corresponds a global section $\nu \in H^0(\SN_{\varphi})$.
Let $\nu_1=\Phi_1 \nu$ and $\nu_2=\Phi_2 \nu$, then we have homomorphisms
\begin{equation}\label{nu-1nu-2}
\SI_{Y,Z}/\SI_{Y,Z}^2 \overset{\nu_1}\longrightarrow \SO_Y 
\quad \mathrm{and} \quad
\SI_{Y,Z}/\SI_{Y,Z}^2 \overset{\nu_2}\longrightarrow \SE.
\end{equation}
Denote $\delta(\nu_1,\nu_2)=(\zeta_1,\zeta_2)$.
There is a first--order deformation $\bar{Y} \overset{\bar{i}}\hookrightarrow Z \times \Delta$ of $Y$ in $Z$ associated with $\nu_1$.
The subscheme $\bar{Y}$ of $Z \times \Delta$ is defined by the ideal $\bar{\SI} \subset \SO_Z \oplus \SO_Z \epsilon$ locally given, over any open affine subset $W=\mathrm{Spec}\,A \subset Z$, by
\begin{equation}\label{barI}
\bar{\SI}=\{a+a'\epsilon \, | \, a \in \SI_{Y,Z} \,\,\mathrm{and}\,\, \nu_1'(a)=-i^{\sharp}a'\},
\end{equation}
where $\nu_1'$ is the composition $\SI_{Y,Z}\to \SI_{Y,Z}/\SI_{Y,Z}^2 \overset{\nu_1}\to \SO_Y$ and $\SO_Z \overset{i^{\sharp}}\to i_*\SO_Y$ is the map of $k$--algebras associated with $Y\overset{i}\hookrightarrow Z$.
Let $\bar{Y}\overset{\widehat{i}}\to Z$ be the composition $\bar{Y} \overset{\bar{i}}\hookrightarrow Z \times \Delta \to Z$.
We see that $(\bar{Y},\widehat{i}\,)$ is an extension of $Y \overset{i}\to Z$. 
Thus, from Proposition~\ref{extension}, we have a commutative diagram like~\eqref{tau^omega}
\begin{equation}\label{e1-omega1-i-nu1}
\xymatrix@C-5pt@R-10pt{
0\ar[r]&\SI/\SI^2 \ar[d]_{\tau^{\omega_1}}\ar[r]&i^*\Omega_Z \ar[d]_-{\omega_1} \ar[r]^-{\mathrm{D}i}&\Omega_Y \ar@{=}[d]\ar[r]&0\\
0 \ar[r] & \SO_Y \ar[r]^-{j_1} & \SF_1 \ar[r]^-{p_1} &  \Omega_Y \ar[r] &  0,}
\end{equation}
where the class of the lower exact sequence $e_1$ is the extension class of $\bar{Y}$.
Using for $\bar{i}^{\sharp}$ a local formula like~\eqref{phi-tilde-local}, we see that $\tau^{\omega_1}=\nu_1$ and then $\zeta_1=[e_1]$.\\
From Proposition~\ref{extension}, to $\nu_2$ corresponds a unique pair $(\WY,\wi)$, where $\WY$ is an $n$--rope with conormal bundle $\SE$ and $\WY \overset{\wi}\to Z$ is a morphism extending $Y \overset{i}\hookrightarrow Z$.
The extension class of $\WY$ is represented by the lower exact sequence $e_2$ defined by the push--out diagram
\begin{equation}\label{e2-omega2-i}
\xymatrix@C-5pt@R-10pt{
0\ar[r]&\SI/\SI^2 \ar[d]_{\nu_2}\ar[r]&i^*\Omega_Z \ar[d]_-{\omega_2} \ar[r]^-{\mathrm{D}i}&\Omega_Y \ar@{=}[d]\ar[r]&0\\
0 \ar[r] & \SE \ar[r]^-{j_2} & \SF_2 \ar[r]^-{p_2} &  \Omega_Y \ar[r] &  0.}
\end{equation}
Then we have $\zeta_2 = [e_2]$.\\
Let $[(e,\omega)]$ be the class defined by a commutative diagram like~\eqref{e-omega-phi} which corresponds to $(\WX,\wphi)$.
Then, from Proposition~\ref{squares}, we have $\delta_1 \nu =[e]$.
Therefore, from the commutativity of~\eqref{construction1}, we have $\Psi_1 [e] =[e_1]$ and $\Psi_2 [e] = [e_2]$.
This proves (1).\\
Now we prove (2).
From Lemma~\ref{glue-ropes}, we know that $\bar{Y}\underset{\scriptscriptstyle{Y}}{\cup}\WY$ is an $(n+1)$--rope with conormal $\SO_Y \oplus \SE$, having a nontrivial structure of $\Delta$--scheme, and there is a unique morphism $\bar{Y}\underset{\scriptscriptstyle{Y}}{\cup}\WY \overset{\widehat{\iota}}\to Z$ extending both $\bar{Y} \overset{\widehat{i}}\to Z$ and $\WY \overset{\wi}\to Z$.
Therefore there is a $\Delta$--morphism $\bar{Y}\underset{\scriptscriptstyle{Y}}{\cup}\WY \overset{\widetilde{\iota}}\to Z \times \Delta$.
Arguing with formulae like~\eqref{local-wphi-whphi}, it is easy to verify that $\bar{Y}\underset{\scriptscriptstyle{Y}}{\cup}\WY \overset{\widetilde{\iota}}\to Z \times \Delta$ is an extension of both $\bar{Y} \overset{\bar{i}}\hookrightarrow Z \times \Delta$ and $\WY \overset{\wi}\to Z$. 
The uniqueness of $\bar{Y}\underset{\scriptscriptstyle{Y}}{\cup}\WY \overset{\widetilde{\iota}}\to Z \times \Delta$ follows from the uniqueness of $\bar{Y}\underset{\scriptscriptstyle{Y}}{\cup}\WY \overset{\widehat{\iota}}\to Z$.
Furthermore, $\bar{Y}\underset{\scriptscriptstyle{Y}}{\cup}\,\WY \overset{\widehat{\iota}}\to Z\,$ is the morphism extending $Y \overset{i}{\to} Z\,$ that, from Proposition~\ref{extension}, corresponds to $\SI_{Y,Z}/\SI_{Y,Z}^2 \overset{\nu_1 \oplus \nu_2}\longrightarrow \SO_Y\oplus \SE$.
Indeed, from~\eqref{e1-omega1-i-nu1}~and~\eqref{e2-omega2-i} there is a commutative diagram
\begin{equation}\label{nus-diag-F}
\xymatrix@C-5pt@R-10pt{
0\ar[r]&\SI/\SI^2 \ar[d]_{\nu_1\oplus\nu_2}\ar[r]&i^*\Omega_Z \ar[d]_-{\omega_3} \ar[r]^-{\mathrm{D}i}&\Omega_Y \ar@{=}[d]\ar[r]&0\\
0 \ar[r] & \SO_Y \oplus \SE \ar[r] & \SF \ar[r] &  \Omega_Y \ar[r] &  0,}
\end{equation}
where the lower exact sequence is the extension class of $\bar{Y}\underset{\scriptscriptstyle{Y}}{\cup}\WY$. 
Therefore, from Proposition~\ref{extension}, there is a morphism $\bar{Y}\underset{\scriptscriptstyle{Y}}{\cup}\WY \overset{\widehat{\iota}\,'}\to Z\,$ extending $Y \overset{i}{\to} Z\,$ defined by $\nu_1\oplus\nu_2$.
Pushing~\eqref{nus-diag-F} out by the projections from $\SO_Y \oplus \SE$ to the factors, we recover diagrams~\eqref{e1-omega1-i-nu1}~and~\eqref{e2-omega2-i}.
So $\widehat{\iota}\,'$ restricts to $\bar{Y} \overset{\widehat{i}}\to Z$ and $\WY \overset{\wi}\to Z$ and therefore we have $\widehat{\iota}\,'=\widehat{\iota}$.\\
Now we prove (3).
We can write the algebra of $\bar{Y}\underset{\scriptscriptstyle{Y}}{\cup}\WY$, over an open affine subset $U=W \cap Y$, as
\begin{equation*}\label{O-U}
\SO=\{(c,f_1,f_2) \in \SO_Y \oplus \SF_1 \oplus \SF_2 \; | \; \mathrm{d}c =p_1f_1=p_2f_2\}.
\end{equation*}
Now, with formulae like~\eqref{local-wphi-whphi} and~\eqref{wi^sharp-formula}
we see that $\widetilde{\iota}\,$ is locally given by
\begin{equation*}\label{widetilde-iota-sharp}
\begin{aligned}
\SO_Z \oplus \SO_Z\epsilon & \overset{\widetilde{\iota}^{\sharp}}{\longrightarrow} i_*\SO\\
\widetilde{\iota}^{\sharp}(a+a'\epsilon)& =(i^{\sharp}a,\, \omega_1(\mathrm{d}a\otimes 1) +j_1 i^{\sharp}a', \, \omega_2(\mathrm{d}a\otimes 1)).
\end{aligned}
\end{equation*}
So we have
\begin{equation}\label{ker-wide-iota}
\mathrm{ker}\, \widetilde{\iota}^{\sharp}=\{a+a'\,\epsilon \;|\; a\in \SI, \,\nu_1(\bar{a})=-i^{\sharp}a' \; \mathrm{and} \; \nu_2(\bar{a})=0\}.
\end{equation}
Let $\SJ$ denote the ideal in $Z$ of the image of the morphism $\WY \overset{\wi}\to Z$.
We know, from Proposition~\ref{extension}, that $\SJ$ is the kernel of $\SI_{Y,Z} \to \SI_{Y,Z}/\SI_{Y,Z}^2 \overset{\nu_2}\to \SE$.
So, from~\eqref{barI}~and~\eqref{ker-wide-iota}, we see that
\begin{equation*}\label{ker-wide-iota2}
\mathrm{ker}\, \widetilde{\iota}^{\sharp}=\bar{\SI} \cap (\SJ +\SO_Z \,\epsilon).
\end{equation*}
This proves (3).\\
Now we define the morphism $\WX \overset{\widetilde{\pi}}\to \bar{Y}\underset{\scriptscriptstyle{Y}}{\cup}\WY$.
Let $\pi'$ denote $X \overset{\pi}\to Y \hookrightarrow \bar{Y}\underset{\scriptscriptstyle{Y}}{\cup}\WY$.
Then we have $\pi'^*(\Omega_{\bar{Y}\underset{\scriptscriptstyle{Y}}{\cup}\WY}) \simeq \pi^*\SF$.
Therefore, from~\ref{noname1.4}, the morphisms $\WX \overset{\widetilde{\pi}}\to \bar{Y}\underset{\scriptscriptstyle{Y}}{\cup}\WY$ extending $X \overset{\pi'}\longrightarrow \bar{Y}\underset{\scriptscriptstyle{Y}}{\cup}\WY$ are in bijection with the commutative diagrams
\begin{equation}\label{e-omega'1}
\xymatrix@C-5pt@R-9pt{
&&&\pi^*\SF \ar[dl]_-{\omega'} \ar[d]^-{\mathrm{D}\pi'}&\\
0 \ar[r] & \SO_X \ar[r] & \SG \ar[r] &  \Omega_X \ar[r] &  0.}
\end{equation}
We construct a diagram like~\eqref{e-omega'1}.
Consider the homomorphism $\SI_{Y,Z}/\SI_{Y,Z}^2 \overset{\nu_1 \oplus \nu_2}\longrightarrow \pi_*\SO_X$.
With the notations of~\eqref{deltas-square} we have $\alpha^{-1}(\nu_1\oplus \nu_2) =\phi(\nu)$.
Hence $\phi(\nu)$ is the composition
\begin{equation*}\label{phi-nu}
\xymatrix@C+10pt{
\pi^*\SI/\SI^2 \ar[r]^-{\pi^*(\nu_1\oplus\nu_2)}& \pi^*\pi_*\SO_X \ar[r]& \SO_X.}
\end{equation*}
Now, from~\eqref{deltas-square}, we obtain $\beta^{-1}\delta_3 (\nu_1\oplus\nu_2) = \delta_2 \phi(\nu)$.
Therefore we obtain a commutative diagram
\begin{equation}\label{phi-SH}
\xymatrix@C-5pt@R-10pt{
0\ar[r]&\pi^*\SI/\SI^2 \ar[d]_{\phi(\nu)}\ar[r]& \varphi^*\Omega_Z \ar[d] \ar[r]^-{\pi^*\mathrm{D}i}& \pi^*\Omega_Y \ar@{=}[d]\ar[r]&0\\
0 \ar[r] & \SO_X \ar[r] & \SH \ar[r] &  \pi^*\Omega_Y \ar[r] &  0,}
\end{equation}
where the lower exact sequence represents the class $\beta^{-1}\delta_3(\nu_1\oplus\nu_2)$ and $\varphi^*\Omega_Z \to \SH$ is the composition $\varphi^*\Omega_Z \overset{\pi^*\omega_3}\longrightarrow  \pi^*\SF \to \SH$.
Also from~\eqref{deltas-square}
we see that $\mathrm{d}\pi ([e])$ is the class of the lower exact sequence in~\eqref{phi-SH}.
Therefore we obtain a commutative diagram
\begin{equation*}\label{pi*SF-omega'}
\xymatrix@C-5pt@R-9pt{
0 \ar[r] & \pi^*\pi_*\SO_X \ar[d] \ar[r] & \pi^*\SF \ar[d]_-{\omega'} \ar[r] & \pi^*\Omega_Y \ar[d]^-{\mathrm{D}\pi} \ar[r]& 0\\
0 \ar[r] & \SO_X \ar[r] & \SG \ar[r] & \Omega_X \ar[r] &0.}
\end{equation*}
Now, from the definition of $\pi'$, we see that the composition $\pi^*\SF \to \pi^*\Omega_Y \overset{\mathrm{D}\pi}\longrightarrow \Omega_X$ is $\pi^*\SF \overset{\mathrm{D}\pi'}\longrightarrow \Omega_X$.
This proves the existence of $\WX \overset{\widetilde{\pi}}\to \bar{Y}\underset{\scriptscriptstyle{Y}}{\cup}\WY$ extending $X \overset{\pi'}\longrightarrow \bar{Y}\underset{\scriptscriptstyle{Y}}{\cup}\WY$.
The rest of (2) follows easily.
\end{Proof}\\     
The following result, with the notations of Proposition~\ref{construction} and~\eqref{nu-1nu-2}, describes the image of the map $\widetilde{\varphi}$.
\begin{theorem}\label{image}
Let $X \overset{\varphi}\to Z$ be a morphism from an integral Cohen--Macaulay variety $X$ to a smooth irreducible variety $Z$.
Let $Y$ be its scheme-theoretic image.
Assume that $Y$ is smooth and that $\varphi$ induces a finite morphism $X \overset{\pi}\to Y$.
Let $(\WX,\widetilde{\varphi})$ be a first--order, locally trivial, infinitesimal deformation of $X \overset{\varphi}\to Z$ defined by a global section $\nu$ of $\, \SN_{\varphi}$.
Then,
\begin{enumerate}
\item the central fiber of the image of the morphism $\, \widetilde{\varphi}$ contains $Y$ and is contained in the first infinitesimal neighborhood of $\, Y$ in $Z$ and is equal to the image of the morphism $\WY \to Z$ obtained from $\nu_2$.
More precisely, the ideal of both the central fiber of the image of $\wphi$ and the image of $\, \WY \to Z \,$ is the kernel of the composite homomorphism
$$\SI_{Y,Z} \to \SI_{Y,Z}/\SI_{Y,Z}^2 \overset{\nu_2}\to \SE.$$
\item The image of $\, \widetilde{\varphi}$ is the scheme-theoretic union of its central fiber and the flat deformation of $\, Y$ defined by $\nu_1$ and it is equal to the image of the $\Delta$--morphism $\bar{Y}\underset{\scriptscriptstyle{Y}}{\cup}\WY \overset{\widetilde{\iota}}\to Z \times \Delta$.
\end{enumerate}
\end{theorem} 
\begin{Proof}
Let $\widetilde{\SJ}$ denote the ideal in $Z \times \Delta$ of the image of $\wphi$ and let $\SJ$ denote the ideal in $Z$ of the central fiber $(\mathrm{im}\, \wphi)_0$.
Let $\bar{\SI}$ denote the ideal of $\bar{Y}$ in $Z \times \Delta$.\\
Observe that a locally split short exact sequence of coherent sheaves on a noetherian separated scheme is split over every open affine subset.
So we can cover $Y$ with open affine subsets $U=W \cap Y=\mathrm{Spec}\,A/I$, where $W=\mathrm{Spec}\,A$ are open affine subsets covering $Z$, so that the first--order, locally trivial, deformation  $\WX$ is trivial over the open affine subsets $V=\pi^{-1}U=\mathrm{Spec}\,B$ covering $X$.
Then we have $B=A/I \oplus M$ with $M=\Gamma(U,\mathcal{E})$.
Let $A \overset{\varphi^{\sharp}}{\to} B$ be the ring homomorphism induced by the morphism $\varphi$.
From Proposition~\ref{extension2} and Lemma~\ref{key-map}, to $(\WX, \wphi)$ corresponds a global section $\nu \in H^0(\SN_{\varphi})$ and, from the isomorphism in~\eqref{g-restricted}, we consider $\nu$ as a class $[(g,\rho)] \in D(X, \varphi)$, for a pair $(g,\rho) \in \mathcal{C}^0(\mathcal{V},\mathcal{H}om_{\SO_X}(\varphi^* \Omega_Z,\SO_X)) \times \mathcal{Z}^1(\mathcal{V},\mathcal{H}om_{\SO_X}(\Omega_X,\SO_X))$ such that $\delta g = \rho \, \mathrm{D}\varphi$, where $\mathcal{V}$ is the con\-si\-dered open affine cover on $X$.\\
Then, from~\eqref{phi-tilde-local}, we see that the ring homomorphism induced by the morphism $\wphi$ is written as
\begin{equation}\label{phi-tilde-local-split}
\begin{aligned}
A \oplus A\epsilon & \overset{\wphi^{\sharp}}{\longrightarrow} B \oplus B \epsilon\\
a+a'\epsilon & \mapsto \varphi^{\sharp}a + 
({g}_{{}_{\smV}}(\mathrm{d}a \otimes 1)+\varphi^{\sharp}a')\,\epsilon,
\end{aligned}
\end{equation}
where $\Omega_A \otimes B \overset{{g}_{{}_{\smV}}}\longrightarrow B$ is the homomorphism given by the cochain $g$ over the open set $V \in \SV$.\\
We prove part (1).
The fact that the ideal of the image of the morphism $\WY \to Z$ obtained from $\nu_2$ is the kernel of $\SI_{Y,Z} \to \SI_{Y,Z}/\SI_{Y,Z}^2 \overset{\nu_2}\to \SE$ follows from Proposition~\ref{extension}.\\
Now we prove that $\SJ$ is the kernel of $\SI_{Y,Z} \to \SI_{Y,Z}/\SI_{Y,Z}^2 \overset{\nu_2}\to \SE$.
This can be locally checked. 
Let $I,J \subset A$, respectively, denote the ideals of $Y$ and $\,(\mathrm{im}\,\widetilde{\varphi})_0$ and let $\bar{I}, \WJ \subset A\oplus A\epsilon$, respectively, denote the ideals of $\bar{Y}$ and $\mathrm{im}\,\widetilde{\varphi}$.
By definition $\WJ$ is the kernel of the homomorphism~\eqref{phi-tilde-local-split}.
Therefore
\begin{equation*}\label{jotatilde-0}
\WJ=\{a+a'\epsilon \;\, | \;\, a \in I \; \mathrm{and} \; {g}_{{}_{\smV}}(\mathrm{d}a \otimes 1)+\varphi^{\sharp}a'=0\}.
\end{equation*}
Moreover, taking the central fiber amounts to tensor the rings with $k[\epsilon]/\epsilon k[\epsilon]$, so $J$ is the image of $\WJ \otimes k[\epsilon]/\epsilon k[\epsilon] \to A$.
Hence
\begin{equation*}\label{jota}
J=\{a \in I \;\, | \;\, {g}_{{}_{\smV}}(\mathrm{d}a \otimes 1) \in \; \mathrm{im}\,\varphi^{\sharp}=A/I \subset B \}.
\end{equation*}
From Lemma~\ref{key-map}, we obtain commutative diagrams
\begin{equation}\label{keyfact}
\xymatrix{
I/I^2 \otimes B  \ar[d]_-{{\phi(\nu)|}_V} \ar[r] & \Omega_A \otimes B  \ar[dl]^-{{g}_{{}_{\smV}}}\\
B \, ,&}\qquad
\xymatrix{
I/I^2 \ar[d]_-{{\alpha \phi(\nu)|}_U} \ar[r] & \Omega_A \otimes A/I  \ar[dl]^-{{f}_{{}_{\smU}}} \\
B \, .&}
\end{equation}
From the splitting $B = A/I \oplus M$, we consider ${f}_{{}_{\smU}}$ as a pair $(({f}_{{}_{\smU}})_{{}_1},({f}_{{}_{\smU}})_{{}_2})\in \mathrm{Hom}_{A/I}(\Omega_A \otimes A/I, A/I)\oplus \mathrm{Hom}_{A/I}(\Omega_A \otimes A/I, M)$.
Hence, with the notations of~\eqref{construction1} and~\eqref{nu-1nu-2}, the right--hand side diagram of~\eqref{keyfact} amounts to commutative diagrams
\begin{equation}\label{keyfact-2}
\xymatrix{
I/I^2 \ar[d]_-{\nu_1} \ar[r] & \Omega_A \otimes A/I  \ar[dl]^-{({f}_{{}_{\smU}})_{{}_1}} \\
A/I \, , & }
\qquad \qquad
\xymatrix{
I/I^2 \ar[d]_-{\nu_2} \ar[r] & \Omega_A \otimes A/I  \ar[dl]^-{({f}_{{}_{\smU}})_{{}_2}} \\
M \, . &}
\end{equation}
From~\eqref{keyfact}~and~\eqref{keyfact-2} we can rewrite $\WJ$ as
\begin{equation}\label{jotatilde}
\WJ=\{a+a'\epsilon \, | \; a \in I, \nu_1(\bar{a})=-\bar{a'} \; \mathrm{and} \; \nu_2(\bar{a})=0\}.
\end{equation}
From~\eqref{jotatilde} and the fact that $J$ is the image of $\WJ \otimes k[\epsilon]/\epsilon k[\epsilon] \to A$ we see that
\begin{equation}\label{jota-2}
J=\{a \in I \; | \; \nu_2(\bar{a})=0\}.
\end{equation}  
Now we prove part (2).
The fact that the image of $\widetilde{\varphi}$ is the scheme-theoretic union of its central fiber and the flat deformation of $Y$ defined by $\nu_1$ is the identity $\widetilde{\SJ}=\bar{\SI}\cap (\SJ +\SO_Z \, \epsilon)$, or locally
\begin{equation}\label{JIcapJA}
\widetilde{J}=\bar{I}\cap (J+A\epsilon).
\end{equation}
We prove~\eqref{JIcapJA}.
The flat subscheme $\bar{Y} \subset Z \times \Delta$ is the image of the deformation of the inclusion morphism $Y \hookrightarrow Z$ defined, from Proposition~\ref{extension2}, by the global section $\nu_1 \in H^0(\SN_{Y,Z})$.
Hence
\begin{equation}\label{itilde}
\bar{I}=\{a+a'\epsilon \, | \; a \in I \; \mathrm{and} \; {({f}_{{}_{\smU}})_{{}_1}}(\mathrm{d}a \otimes 1)=-{\bar a'}\} =\{a+a'\epsilon \, | \; a \in I \; \mathrm{and} \; \nu_1(\bar{a})=-{\bar a}'\}.
\end{equation}
From~\eqref{jota-2},~\eqref{itilde} and~\eqref{jotatilde} we obtain~\eqref{JIcapJA}.\\
The fact that the image of the $\Delta$--morphism $\bar{Y}\underset{\scriptscriptstyle{Y}}{\cup}\WY \overset{\widetilde{\iota}}\to Z \times \Delta$ is the scheme-theoretic union of the image of $\WY \to Z$ and $\bar{Y}$ is proven in Proposition~\ref{construction}.
\end{Proof}\\
Now we state the main result of this section.
\begin{theorem}\label{main}
Let $X \overset{\varphi}\to Z$ be a morphism from an integral Cohen--Macaulay variety $X$ to a smooth irreducible variety $Z$.
Let $Y$ be its scheme-theoretic image.
Assume that $Y$ is smooth and that $\varphi$ induces a finite morphism $X \overset{\pi}\to Y$.
Let $\SE$ be the locally free $\SO_Y$--module $\pi_*\SO_X/\SO_Y$.
If the map $H^1(\SN_{\pi}) \to H^1(\SN_{\varphi})$ is injective then every rope over $Y$, with conormal bundle contained in $\SE$, embedded in $Z$ is the central fiber of the image of some first--order, locally trivial, infinitesimal deformation of $\varphi$.
\end{theorem}
\begin{Proof}
From~\ref{noname1.3}, a rope $\WY$ over $Y$ with conormal bundle $\SE'$, embedded in $Z$ corresponds to a surjective map $\SI_{Y,Z}/\SI_{Y,Z}^2 \overset{\tau}\twoheadrightarrow \SE'$ and its ideal $\SJ$ is the kernel of $\SI \twoheadrightarrow \SI_{Y,Z}/\SI_{Y,Z}^2 \overset{\tau}\twoheadrightarrow \SE'$.\\
Assume that $\SE' \subset \SE$ and let $\nu_2$ denote the induced map $\SI_{Y,Z}/\SI_{Y,Z}^2 \overset{\nu_2}\to \SE$.
The hypothesis that $H^1(\SN_{\pi}) \to H^1(\SN_{\varphi})$ is injective is equivalent, from the cohomology sequence of~\eqref{normal-extension}, to the fact that the map
\begin{equation*}\label{lifting}
H^0(\SN_{\varphi}) \overset{\phi}\longrightarrow H^0(\pi^* \SN_{Y,Z})
\end{equation*}
in~\eqref{deltas-square} is surjective.
Therefore the element $(0,\nu_2) \in H^0(\pi^*\SN_{Y,Z})$ admits a lifting to a global section $\nu \in H^0(\SN_{\varphi})$.
Let $\widetilde{\varphi}$ be the associated first order, locally trivial, infinitesimal deformation.
From Theorem~\ref{image}, we see that the central fiber $(\mathrm{im}\,\widetilde{\varphi})_0$ has ideal $\SJ$.
\end{Proof}\\
As a consequence we obtain the following result, from which we obtain Theorem~\ref{key.curves.0} which is the key infinitesimal result we will use in the proof of the embedded smoothing Theorem~\ref{embsmoothing} of ribbons over curves.
\begin{theorem}\label{key-curves}
Let $X \overset{\varphi}\to Z$ be a morphism from a smooth irreducible curve $X$ to a smooth irreducible variety $Z$.
Let $Y$ be its scheme-theoretic image.
Assume that $Y$ is smooth and that $\varphi$ induces a finite morphism $X \overset{\pi}\to Y$.
Let $\SE$ be the locally free $\SO_Y$--module $\pi_*\SO_X/\SO_Y$.
Then every rope over the smooth irreducible curve $Y$, with conormal bundle contained in $\SE$, embedded in $Z$ is the central fiber of the image of some first--order infinitesimal deformation of $\varphi$.
\end{theorem}
\begin{Proof}
For $X$ is a curve we see that the support of $\SN_{\pi}$ is a finite set.
Therefore $H^1(\SN_{\pi})=0$.
\end{Proof}\\
Finally, keeping previous notations, we identify when we obtain the subvariety $Y$ itself as central fiber of the image of $\wphi \,$ and when $\wphi$ factors through $Y \times \Delta$.
\begin{proposition}\label{equivalence}
In the situation of Theorem~\ref{image}, the following are equivalent conditions:\\
$\mathrm{(1)}\,\, Y=(\mathrm{im}\, \widetilde{\varphi})_0$.
$\mathrm{(2)}\,\, \nu_2=0$.
$\mathrm{(2')}\,\, (\WY, \wi)$ is the pair consisting of the split rope and its projection to $Y$.
$\mathrm{(3)}\,\, \mathrm{im}\, \widetilde{\varphi}\subset \bar{Y}$.
$\mathrm{(4)}\,\, \mathrm{im}\, \widetilde{\varphi}=\bar{Y}$.
$\mathrm{(5)}\,\, \mathrm{im}\, \widetilde{\varphi}$ is flat over $\Delta$.
\end{proposition}
\begin{Proof} Apart from (5) $\Rightarrow$ (1), all equivalences are direct consequence of Proposition~\ref{extension}, Theorem~\ref{image}, \eqref{jotatilde}~and~\eqref{itilde}.\\
(5) $\Rightarrow$ (1) From the ``Local criteria of flatness" \cite[20.C]{Matsumura80}, we see that $\mathrm{im}\,\widetilde{\varphi}$ is flat over $\Delta$ iff the surjective map $\widetilde{\SJ} \otimes k[\epsilon]/\epsilon k[\epsilon] \to \SJ$ is an isomorphism.
Moreover, it is easy to see that if $\widetilde{\SJ} \otimes k[\epsilon]/\epsilon k[\epsilon] \to \SJ$ is injective then $\widetilde{\SJ} \cap \SO_{Z} \epsilon \subset \SJ \epsilon$.
Now, from~\eqref{jotatilde}, we obtain $\SI \epsilon = \widetilde{\SJ} \cap \SO_{Z} \epsilon$.
So we have $\SI \epsilon \subset \SJ \epsilon$ and therefore $\SI \subset \SJ$.
\end{Proof}
\begin{proposition}
In the situation of Theorem~\ref{image}, the following are equivalent conditions:\\
$\mathrm{(1)}\,\, (\WX, \widetilde{\varphi})$ factors through a first--order, locally trivial, infinitesimal deformation of $\pi$.
$\mathrm{(2)}\,\, \nu \in H^0(\SN_{\pi})$.
$\mathrm{(3)}\,\, \bar{Y}= Y \times \Delta$ and $(\WY, \wi)$ is the pair consisting of the split rope and its projection to $Y$.{\hfill $\square$}
\end{proposition}
\section{Embedding in projective space for ribbons over curves}
From now on, $Y$ will be a smooth irreducible projective curve of arbitrary genus $g$, and $\SE$ a line bundle on $Y$.
We consider ribbons over $Y$ with conormal bundle $\SE$.\\
Recall, from Proposition~\ref{extension}, that maps from ribbons with conormal bundle $\SE$, extending a fixed closed immersion $Y \overset{i}\hookrightarrow \PP^r$, are in correspondence with homomorphisms $\SN_{Y,\PP^r}^* \to \SE$ so that closed immersions correspond to the surjective ones.\\
In the next two results we obtain a criterion to decide whether all ribbons with fixed conormal bundle $\SE$ can be embedded in the same projective space with support in a fixed embedding of the reduced part $Y$.
\begin{lemma}\label{epiisopen}
Let $Y$ be a smooth irreducible projective curve, let $\SE$ be a line bundle on $Y$ and let $\SF$ be a locally free sheaf of finite rank on $Y$.
\begin{enumerate}
\item The surjective homomorphisms from $\SF$ to $\SE$ form an open set of 
$\, \mathrm{Hom}(\SF,\SE)$ which is the complement of an algebraic cone.
\item Consider an extension of a coherent sheaf $\SF'$ by $\SF$.
Let $\, \mathrm{Hom}(\SF,\SE) \overset{\delta}\to \mathrm{Ext}^1(\SF',\SE)$ be the connecting map induced by the extension.
If the class of the split extension lifts to an epimorphism, then every class in the image of $\,\delta$ can be lifted to an epimorphism.
\end{enumerate}
\end{lemma}
\begin{Proof}
(1) A homomorphism $\SF \overset{\phi}\to \SE$ is not surjective at a point $y \in Y$ iff the induced homomorphism $\SF \otimes k(y) \overset{\phi(y)}\longrightarrow \SE \otimes k(y)$ vanishes.\\
Let consider the closed subscheme $\,\Gamma =\{(y,\phi) \,\, | \,\, \phi(y)=0\} \subset Y \times \mathrm{Hom}(\SF,\SE)$.
The projection $B$ of $\, \Gamma$ to $\mathrm{Hom}(\SF,\SE)$ will then be a closed set and it is clearly a cone.\\
(2) Assume that for a class $\zeta \in \mathrm{Ext}^1 (\SF',\SE)-\{0\}$ in the image of $\delta$, every lifting of $\zeta$ by $\delta$ lies in the closed algebraic cone $B$ formed by the non-surjective homomorphisms.
Fix $s \in B$ with $\delta(s)= \zeta$.
Then for every $\lambda \in k^*$ and every $v \in \mathrm{ker}\,\delta$ the identity $\delta (s+\lambda^{-1}v)= \zeta $ and our assumption imply that $s+\lambda^{-1}v \in B$. Since $B$ is a cone, we deduce that for every $\lambda \in k^*$, $B$ contains the set $\lambda s+ \mathrm{ker}\,\delta$.
Since $B$ is closed, we deduce that $B$ also contains the kernel of $\delta$.
So we get a contradiction.
\end{Proof}
\begin{proposition}\label{nondegembedding}
Let $Y$ be a smooth irreducible projective curve in $\PP^r$, with $r \geq 3$, and let $\SE$ be a line bundle on $Y$.
\begin{enumerate}
\item If the split ribbon with conormal bundle $\SE$ can be embedded as a nondegenerate closed subscheme in $\PP^r$ extending the embedding $\, Y \hookrightarrow \PP^r$, then every ribbon with conormal bundle $\SE$ that maps to $\PP^r$ extending $\, Y \hookrightarrow \PP^r$ can be embedded as a nondegenerate closed subscheme in $\PP^r$ extending $\, Y \hookrightarrow \PP^r$.
\item If the connecting map $\mathrm{Hom}(\SN_{Y,\PP^r}^*, \SE) \overset{\delta}\to \mathrm{Ext}^1(\Omega_Y,\SE)$ is surjective and the split ribbon with conormal bundle $\SE$ can be embedded as a nondegenerate closed subscheme in $\PP^r$ extending $\, Y \hookrightarrow \PP^r$, then every ribbon with conormal bundle $\SE$ can be embedded as a nondegenerate closed subscheme in $\PP^r$ extending $\, Y \hookrightarrow \PP^r$.
\end{enumerate}
\end{proposition}
\begin{Proof}
If $\mathrm{Hom}(\SN_{Y,\PP^r}^*, \SE) \overset{\delta}\to \mathrm{Ext}^1(\Omega_Y,\SE)$ is surjective then, from Proposition~\ref{extension}, every ribbon with conormal bundle $\SE$ is mapped to $\PP^r$ extending $\, Y \hookrightarrow \PP^r$.
So the second part follows from the first.\\
Let $L$ be a hyperplane in $\PP^r$.
If $Y$ is contained in $L$ then $\mathrm{Hom}(\SN_{Y, L}^*, \SE)$ is a subspace of $\mathrm{Hom}(\SN_{Y, \PP^r}^*, \SE)$.
Take $s \in \mathrm{Hom}(\SN_{Y, L}^*, \SE)$.
From~\eqref{wi^sharp-formula}~and~\eqref{push-tau}, we easily see that the morphism $\WY \to \PP^r$ defined by $s$ is the composition of the morphism $\WY \to L$ defined by $s$ as an element of $\mathrm{Hom}(\SN_{Y, L}^*, \SE)$ and $L \hookrightarrow \PP^r$.
Therefore the image of the extension morphism $\WY \to \PP^r$ defined, from Proposition~\ref{extension}, by an element $s \in \mathrm{Hom}(\SN_{Y, \PP^r}^*, \SE)$ is contained in a hyperplane $L$ iff $Y$ is contained in $L$ and $s \in \mathrm{Hom}(\SN_{Y, L}^*, \SE)$.
So a nondegenerate embedding in $\PP^r$ of the ribbon asso\-ciated to an element $\zeta \in \mathrm{Ext}^1(\Omega_Y,\SE)$ corresponds to a lifting of $\zeta$ to a surjective homomorphism $s \in \mathrm{Hom}(\SN_{Y, \PP^r}^*, \SE)$ that does not belong to any subspace $\mathrm{Hom}(\SN_{Y, L}^*, \SE)$ when $Y \subset L$ and $L$ is a hyperplane.
Our hypothesis says that this lifting exists for $\zeta=0$.\\
The union of the subspaces $\mathrm{Hom}(\SN_{Y, L}^*, \SE)$, where $L$ is any hyperplane in $\PP^r$ containing $Y$, is a closed algebraic cone.
From Lemma~\ref{epiisopen}, the set of non--surjective homomorphisms is a closed algebraic cone.
So the union of both cones is a closed algebraic cone $B'$ in $\mathrm{Hom}(\SN_{Y, \PP^r}^*, \SE)$.
Now, arguing like in the proof of Lemma~\ref{epiisopen}~(2), we see that every class $\zeta \in \mathrm{Ext}^1(\Omega_Y,\SE)$ in the image of $\delta$ can be lifted to an element of $\mathrm{Hom}(\SN_{Y, \PP^r}^*, \SE)$ that does not belong to $B'$.
This element defines a nondegenerate embedding in $\PP^r$, extending $\, Y \hookrightarrow \PP^r$, of the ribbon with extension class $\zeta$.
\end{Proof}
\begin{remark}\label{deltasurjectivity}
{\rm The map $\delta$ is surjective if $H^1(\SE \otimes \SO_Y(1))=0$.}{\hfill $\square$}
\end{remark}
Our next goal will be to find nondegenerate projective embeddings, in the same projective space, for all ribbons with conormal bundle $\SE$ supported over a (possibly degenerate) projective embedding of the base curve $Y$.
According to Proposition~\ref{nondegembedding}, we first look for a nondegenerate projective embedding of the split ribbon with conormal bundle $\SE$.
The method we will use is suggested by the following proposition.
\begin{proposition}\label{ruled}
Let $Y$ be a smooth irreducible projective curve and let $Y \overset{i}\hookrightarrow \PP^r$ be a closed immersion, with $r \geq 3$.
Let $\WY$ be the split ribbon over $Y$ with conormal bundle $\SE$.
Assume that there is a closed immersion $\WY \overset{\wi}\hookrightarrow \PP^r$ extending $Y \overset{i}\hookrightarrow \PP^r$.
Then
\begin{enumerate}
\item there is an extension $\;\, 0 \to \SE \otimes \SO_Y(1) \to \SM \to  \SO_Y(1) \to 0 \;\,$ such that the split ribbon $\WY$ is the first infinitesimal neighborhood of the section defined by the surjective map $\SM \to \SO_Y(1)$ inside the geometrically ruled surface $S=\PP_Y(\SM)$.
\item Let $\SO_S(1)$ denote the fundamental line bundle on $S=\PP_Y(\SM)$.
There is a morphism $S \overset{\psi}\to \PP^r$, with $\psi^*\SO_{\PP^r}(1)=\SO_S(1)$, whose image is a, possibly singular, scroll and whose restriction to $\,\WY$ is the embedding $\WY \overset{\wi}\hookrightarrow \PP^r$.
Moreover, we have $H^0(\SO_S(1))= H^0({\SO_S(1) |}_{\WY})$.
\end{enumerate}
\end{proposition}
\begin{Proof}
The ribbon embedded in $\PP^r$ by a surjection $\SN_{Y,\PP^r}^* \overset{\tau}\to \SE$ is the split ribbon iff this surjection extends to a (surjective) map $i^*\Omega_{\PP^r}\overset{\tau}\to \SE$.\\
Now assume that the last surjection exists.
Let $\SF$ be the kernel of $i^*\Omega_{\PP^r}(1) \overset{\tau \otimes 1}\longrightarrow \SE \otimes \SO_Y(1)$.
We consider the pullback to $Y$ of the Euler sequence in $\PP^r$.
Let $\SM$ denote the cokernel of the composition $\SF \to i^*\Omega_{\PP^r}(1) \to H^0(\SO_{\PP^r}(1))\otimes \SO_Y$.
So we obtain the exact sequence in (1) according to a diagram
\begin{equation*}\label{rank2ext-diagram}
\xymatrix@C-10pt@R-10pt{
& &   0 \ar[d] & 0 \ar[d] &  \\
0 \ar[r] & \SF \ar@{=}[d] \ar[r] &  i^*\Omega_{\PP^r}(1) \ar[d] \ar[r]^-{\tau \otimes 1} &
\SE \otimes \SO_Y(1) \ar[d] \ar[r] & 0 \\
0 \ar[r] & \SF \ar[r] & H^0(\SO_{\PP^r}(1))\otimes\SO_Y \ar[d] \ar[r] & \SM \ar[d] \ar[r]& 0 \\
 &  & \SO_Y(1) \ar[d] \ar@{=}[r] & \SO_Y(1) \ar[d]& \\
 & &  0 & 0 &.
}
\end{equation*}
Let $S=\PP_Y(\SM)$ be the geometrically ruled surface associated with $\SM$ with projection $S \overset{p}\to Y$.
Let $Y \hookrightarrow S$ be the section associated with the quotient $\SM \to \SO_Y(1)$. Then we have (see \cite[V 2.6]{Hartshorne77}) $p^*(\SE \otimes \SO_Y(1))=\SO_S(1) \otimes \SO_S(- Y)$ and hence ${\SO_{S}(-Y) |}_Y=\SE$.
The last identity means that the first infinitesimal neighborhood of $Y$ inside $S$ is a ribbon $\WY$ over $Y$ with conormal bundle $\SE$. It is the split ribbon because the restriction of $p$ to $\WY$ gives a retraction from $\WY$ to $Y$.
This proves (1).\\
Pulling back to $S$ the surjective map $H^0(\SO_{\PP^r}(1))\otimes \SO_Y \to \SM$, and composing with the canonical surjection $p^*\SM \to \SO_S(1)$, we obtain a surjective map $H^0(\SO_{\PP^r}(1))\otimes\SO_S \to \SO_S(1)$ that defines a morphism $S \overset{\psi} \to \PP^r$.
By construction, the composition of $Y \hookrightarrow S$ and $S \overset{\psi} \to \PP^r$ is the given map $Y \overset{i}\hookrightarrow \PP^r$.
Moreover, let $\WY \overset{\wi'} \to \PP^r$ be the restriction to $\,\WY$ of $\,S \overset{\psi}\to \PP^r$.
The map $\wi'$ corresponds to the composition $\SN_{Y,\PP^r}^* \to \SN_{Y,S}^* \overset{\sim}\to \SE$ that, by construction, agrees with the original surjection $\SN_{Y,\PP^r}^* \overset{\tau}\to \SE$.
So $\WY \overset{\wi'} \to \PP^r$ is the embedding $\WY \overset{\wi}\to \PP^r$ and therefore the image of $S$ must be a surface since it is reduced and contains $\WY$.\\
On the other hand, the sheaves $R^ip_*(\SO_S(-2Y)\otimes \SO_S(1))$ vanish for $i=0,1$.
So we have $p_*\SO_S(1) \simeq p_*({\SO_S(1)|}_{\WY})$ and therefore $H^0(\SO_S(1))= H^0({\SO_S(1) |}_{\WY})$.
This ends the proof of (2).
\end{Proof}
\begin{proposition}\label{ruled2}
Let $Y$ be a smooth irreducible projective curve and let $S=\PP_Y(\SN) \overset{p} \to Y$ be the geometrically ruled surface associated with a rank two locally free sheaf $\SN$ over $Y$.
Let $\SO_S(1)$ denote the fundamental line bundle on $S=\PP_Y(\SN)$.
\begin{enumerate}
\item A section $Y \hookrightarrow S$ whose first infinitesimal neighborhood inside $S$ is the split ribbon $\WY$ over $Y$ with conormal bundle $\SE$ is equi\-valent to a line sub--bundle $\SL \hookrightarrow \SN$ such that there is an exact sequence $0 \to \SL \to \SN \to \SE^{-1} \otimes \SL \to 0$.
In this case $\SO_S(Y)=\SO_S(1) \otimes p^*\SL^{-1}$.
\item Let $\SO_Y(1)$ be a very ample line bundle on $Y$.
Let $\, \PP^s \subset \PP^r$ denote, respectively, the projective spaces of (one quotients) of $\,H^0(\SO_Y(1))$ and $H^0(\SO_Y(1)) \oplus H^0(\SE \otimes \SO_Y(1))$.
Let $Y \overset{i}\hookrightarrow \PP^r$ be the closed immersion defined as the composition of the embedding $Y \hookrightarrow \PP^s$ given by the complete linear series $H^0(\SO_Y(1))$ and $\PP^s \subset \PP^r$.\\
Assume that $H^1(\SE \otimes \SO_Y(1))=0$.
Assume that there is an exact sequence
\begin{equation}\label{ruled2-0}
0 \to \SL \to \SN \to \SE^{-1} \otimes \SL \to 0,
\end{equation}
with $\SL$ a line bundle, and let $Y \hookrightarrow S$ be the associated section.\\
If the line bundle $\SO_S(1)\otimes p^*(\SE \otimes \SO_Y(1) \otimes \SL^{-1})$ is globally generated then its complete linear series defines a morphism $S \overset{\psi}\to \PP^r$ such that the composition with the section is the given embedding $Y \overset{i}\hookrightarrow \PP^r$.\\
Moreover, the induced morphism from the split ribbon $\WY$ to $\PP^r$ is defined by the complete li\-near series of ${(\SO_S(1)\otimes p^*(\SE \otimes \SO_Y(1) \otimes \SL^{-1}))|}_{\WY}$ and therefore its image is nondegenerate.
\end{enumerate}
\end{proposition}
\begin{Proof} Let $Y \hookrightarrow S$ be a section defined by a surjective map $\SN \twoheadrightarrow \SK$ with $\SK$ an invertible sheaf.
Let $\SL$ be the kernel of $\SN \twoheadrightarrow \SK$.
Then we have (see \cite[V 2.6]{Hartshorne77}) $p^*\SL = \SO_S(1)\otimes\SO_S(-Y)$ and $\SL =  {(\SO_S(1)\otimes\SO_S(-Y))|}_Y$.
Moreover, if $\WY$ is the first infinitesimal neighborhood of $Y$ inside $S$ then $\WY$ is the split ribbon over $Y$ with conormal bundle ${\SO_S(-Y)|}_Y$.
So we have $\SK =\SE^{-1} \otimes \SL$ iff ${\SO_S(-Y)|}_Y= \SE$.
This proves (1).\\
We prove (2).
Denote $\SL' =  \SE \otimes \SO_Y(1) \otimes \SL^{-1}$.
Observe first that
\begin{equation}\label{ruled2-1}
H^0(\SO_S(1) \otimes p^*\SL')= H^0(\SN \otimes \SL').
\end{equation}
From the definition of the section we have ${\SO_S(1)|}_Y= \SE^{-1} \otimes \SL$.
So we have the identity
\begin{equation*}\label{ruled2-2}
{(\SO_S(1) \otimes p^*\SL')|}_Y= \SO_Y(1)
\end{equation*}
and the exact sequence
\begin{equation}\label{ruled2-2-1}
0 \to \SO_S(-Y)\otimes \SO_S(1)\otimes p^*\SL' \to \SO_S(1)\otimes p^*\SL' \to \SO_Y(1) \to 0.
\end{equation}
Now pushing--down to $Y$ the sequence~\eqref{ruled2-2-1} we obtain
\begin{equation*}\label{ruled2-3}
0 \to \SE \otimes \SO_Y(1) \to \SN \otimes \SL' \to \SO_Y(1) \to 0,
\end{equation*}
which is~\eqref{ruled2-0} twisted by $\SL'$. 
From the assumption $H^1(\SE \otimes \SO_Y(1))=0$ and the isomorphism~\eqref{ruled2-1}, we obtain an exact sequence
\begin{equation*}\label{ruled2-4}
0\to H^0(\SE \otimes \SO_Y(1)) \to H^0(\SO_S(1)\otimes p^*\SL') \to H^0(\SO_Y(1)) \to 0.
\end{equation*}
Therefore we have a commutative diagram
\begin{equation*}\label{ruled2-4-1}
\xymatrix@R-8pt{
H^0(\SO_S(1)\otimes p^*\SL')\otimes \SO_S \ar@{->>}[d] \ar@{->>}[r] & \SO_S(1)\otimes p^*\SL' \ar@{->>}[d]\\
H^0(\SO_Y(1)) \otimes \SO_Y \ar@{->>}[r] & \SO_Y(1)}
\end{equation*}
where the horizontal arrows are the evaluation morphisms.
This proves that there is a morphism $S \overset{\psi} \to \PP^r$ whose restriction to $Y$ is the given embedding $Y \hookrightarrow \PP^r$.\\  
Finally, we have $R^ip_*(\SO_S(-2Y)\otimes \SO_S(1)\otimes p^*\SL')=0$ for $i=0,1$.
So $p_*(\SO_S(1) \otimes p^*\SL')\simeq p_*({(\SO_S(1) \otimes p^*\SL')|}_{\WY})$ and therefore $H^0(\SO_S(1) \otimes p^*\SL') = H^0({(\SO_S(1) \otimes p^*\SL')|}_{\WY})$.
\end{Proof}
\begin{noname}
{\rm For a base curve $Y$ of arbitrary genus $g$, the way to place the split ribbon $\WY$ with conormal bundle $\SE$ is to take the first infinitesimal neighborhood of the section defined by $\, \SO_Y \oplus \SE \twoheadrightarrow \SO_Y$ inside $\PP(\SO_Y \oplus \SE)$.
This is the split ribbon with conormal bundle $\SE$ viewed as (see \cite[1.1]{BayerEisenbud95}) the first infinitesimal neighborhood of the null section in $\mathrm{\bf V}(\SE)$.}
\end{noname}
\begin{theorem}\label{embedding}
Let $Y$ be a smooth irreducible projective curve, let $\SO_Y(1)$ be a very ample line bundle on $Y$ and let $\SE$ be a line bundle on $Y$.
Let $\, \PP^s \subset \PP^r$ denote, respectively, the projective spaces of (one quotients) of $\,H^0(\SO_Y(1))$ and $H^0(\SO_Y(1)) \oplus H^0(\SE \otimes \SO_Y(1))$.
Let $Y \overset{i}\hookrightarrow \PP^r$ be the closed immersion defined as the composition of the embedding $Y \hookrightarrow \PP^s$ given by the complete linear series $H^0(\SO_Y(1))$ and $\PP^s \subset \PP^r$.
Let $d =- \, \mathrm{deg}\, \SE$ and assume that $\mathrm{deg} \, \SO_Y(1) \geq \max \,\{d+2g+1, 2g+1\}$.
Then we have $r \geq 3$ and every ribbon over $Y$ with conormal bundle $\SE$ admits a nondegenerate embedding in $\PP^r$ with degenerate support $Y \overset{i}\hookrightarrow \PP^r$.
\end{theorem}
\begin{Proof}
We can apply Proposition~\ref{ruled2} for the section defined by $\SO_Y \oplus \SE \twoheadrightarrow \SO_Y$ inside $\PP(\SO_Y \oplus \SE)$.
The assumption $\mathrm{deg} \, \SO_Y(1) \geq \max \,\{d+2g+1, 2g+1\}$ implies that the line bundle $\SO_S(1)\otimes p^*\SO_Y(1)$ is very ample (see \cite[V Ex. 2.11]{Hartshorne77}).
So we obtain a nondegenerate embedding in $\PP^r$ for the split ribbon with conormal bundle $\SE$.
Now we obtain Theorem~\ref{embedding} from Remark~\ref{deltasurjectivity} and Proposition~\ref{nondegembedding}.\\
Observe that $h^0(\SO_Y(1)) \geq g+2$ and $h^0(\SE \otimes \SO_Y(1)) \geq g+2$, so $r \geq 3$ and the embedding $Y \overset{i}\hookrightarrow \PP^r$ is degenerate.
\end{Proof}
\begin{remark}\label{elliptic-rational}
{\rm Although with this kind of embeddings we will be able to obtain our main result regarding smoothing of ribbons, in some cases we can obtain nondegenerate embeddings, in the same projective space $\PP^r$, for all ribbons over $Y$ with conormal bundle $\SE$ extending a nondegenerate embedding of $\,Y$ in $\PP^r$.\\
(1) If $Y$ is an elliptic curve and $d =- \, \mathrm{deg}\, \SE \geq 5$ then we can embed all ribbons with conormal bundle $\SE$ in $\PP^{d-1}$ over an elliptic normal curve $Y \subset \PP^{d-1}$ of degree $d$.\\
Indeed, we take a line bundle $\SO_Y(1)$ on $Y$ with $\mathrm{deg}\, \SO_Y(1)=d$ such that $\SO_Y(1) \ncong \SE^{-1}$.
So we have $H^0(\SE \otimes \SO_Y(1))=0$ and $H^1(\SE \otimes \SO_Y(1))=0$.
We consider the embedding $Y \overset{i}\hookrightarrow \PP^{d-1}$ defined by the complete linear series of $\SO_Y(1)$.
To place the split ribbon in $\PP^{d-1}$ over $Y\overset{i}\hookrightarrow \PP^{d-1}$ we use Proposition~\ref{ruled2} looking at surfaces associated with indecomposable rank two vector bundles on $Y$.\\
For $d$ odd ($\geq 5$), we fix a point $O \in Y$ and we consider the non--trivial extension $0 \to \SO_Y \to \SN \to \SO_Y(O)\to 0$.
Let $S=\PP_Y(\SN)$ be the geometrically ruled surface over $Y$ associated with $\SN$.
If we take $\SL$ such that $\SL^2 = \SE(O)$, then $H^0(\SN \otimes \SL^{-1})=H^0(\SL^{-1})\oplus H^0(\SL^{-1}(O))$ and therefore we can take a nowhere vanishing global section of $\SN \otimes \SL^{-1}$ defining an exact sequence $0 \to \SL \to \SN \to \SE^{-1} \otimes \SL \to 0$.
From Proposition~\ref{ruled2}, we get the split ribbon with conormal bundle $\SE$ as the first infinitesimal neighborhood of the section defined by the surjection $\SN \to \SE^{-1} \otimes \SL$ inside $S$.
Moreover, the line bundle $\SO_S(1)\otimes p^*(\SE \otimes \SO_Y(1)\otimes \SL^{-1})$ is very ample (see \cite[V Ex. 2.12]{Hartshorne77}) and its complete linear series gives an embedding for $S$ in $\PP^{d-1}$.
So arguing like in the proof of Theorem~\ref{embedding}, we get our desired embedding for all ribbons with conormal bundle $\SE$ in $\PP^{d-1}$ with nondegenerate embedded reduced support $Y \hookrightarrow \PP^{d-1}$.\\
For $d$ even ($\geq 6$), we consider likewise the geometrically ruled surface $S=\PP_Y(\SN)$ over $Y$ associated with the non--trivial extension $0 \to \SO_Y \to \SN \to \SO_Y \to 0$.
We take $\SL$ such that $\SL^2 = \SE$ and a nowhere vanishing section in $H^0(\SN \otimes \SL^{-1})=H^0(\SL^{-1})\oplus H^0(\SL^{-1})$.\\
For $d=3,4$, again $Y$ elliptic, this does not work but we can place ribbons in $\PP^{d+1}$ over $Y \overset{i}\hookrightarrow \PP^d$.
Now we take $\SO_Y(1)= \SE^{-1}(O)$ and we proceed like above for odd and even cases.\\
For $d=2$, ribbons can be embedded in $\PP^5$ with degree $8$ over $Y \subset \PP^3$.
We take $\SO_Y(1)= \SE^{-1}(2O)$, a line bundle $L \ncong \SO_Y$ with $\mathrm{deg}\, L=0$, and $\SN=\SO_Y \oplus L$.
We also take $\SL$ such that $\SL^2 = \SE \otimes L$ and a nowhere vanishing global section of $\SN \otimes \SL^{-1}$ and proceed like above.\\
(2) If $Y=\PP^1$ and $d \geq 4$ then we can embed all ribbons with conormal bundle $\SO_{\PP^1}(-d)$ in $\PP^{d-1}$ over a rational normal curve $Y$ of degree $d-1$ in $\PP^{d-1}$.\\
Indeed, a nowhere vanishing global section of $\SO_{\PP^1}(d-2) \oplus \SO_{\PP^1}(2)$ gives an exact sequence $0 \to \SO_{\PP^1}(-d+2) \to \SO_{\PP^1}\oplus \SO_{\PP^1}(-d+4) \to \SO_{\PP^1}(2)\to 0$.
Let $S$ be the geometrically ruled surface over $\PP^1$ with invariant $e=d-4$.
The line bundle $\SO_S(1)\otimes p^*\SO_{\PP^1}(d-3)$ is very ample and we can apply Proposition~\ref{ruled2} to embed the split ribbon with conormal $\SO_{\PP^1}(-d)$ inside a rational normal scroll of degree $d-2$ in $\PP^{d-1}$ and with reduced support in a rational normal curve in $\PP^{d-1}$.
Next we use, like above, Remark~\ref{deltasurjectivity} and Proposition~\ref{nondegembedding} to obtain embedding in $\PP^{d-1}$ with nondegenerate embedded reduced support $Y \subset \PP^{d-1}$ for all ribbons with conormal bundle $\SO_{\PP^1}(-d)$.{\hfill $\square$}}
\end{remark}
\section{Smoothing}\label{smoothing}
In this section we show that, under weak conditions on the conormal bundle $\SE$, every ribbon of arithmetic genus greater than or equal to $3$ over a smooth irreducible projective curve of arbitrary genus $g$ (with very few exceptions if $g=0$ or $g=1$) is smoothable.\\
A smoothing of a ribbon $\WY$ over a smooth irreducible projective curve $Y$ is an integral family $\SY$ proper and flat over a smooth pointed affine curve $(T,0)$, whose general fiber $\SY_t, t \neq 0$, is a smooth irreducible projective curve and whose central fiber $\SY_0$ is isomorphic to $\WY$.\\
If $\WY \subset \PP^r$, $\SY \subset \PP_T^r$ as closed subschemes, and $\SY_0=\WY$ we call $\SY$ an embedded smoothing.\\
Our main Theorem~\ref{embsmoothing} gives sufficient conditions for an embedded ribbon $\WY \subset \PP^r$ to have an embedded smoothing.\\
Over a smooth irreducible projective curve $Y$ of arbitrary genus $g$ we consider ribbons with a fixed conormal bundle $\SE$.
Let us denote $d =- \mathrm{deg} \, \SE$. Then the arithmetic genus is $\,p_a(\WY)=d+2g-1$.\\
In Theorem~\ref{embsmoothing} we assume the existence of a smooth irreducible double cover $X \overset{\pi}\to Y$ with $\pi_*\SO_X/\SO_Y = \SE$.
Such a double cover $X$ is determined (see e.g. \cite[I.17]{BPVdV84}) by $\SE$ and its branch locus, an effective divisor, smooth for $X$ to be smooth, with associated line bundle $\SE^{-2}$.
Therefore the existence of such a double cover is equivalent to either
\begin{enumerate}
\item the existence on $Y$ of a non--zero effective reduced divisor with associated line bundle $\SE^{-2}$, or
\item $\SE^{-2}=\SO_Y$ and $H^0(\SE)=0$.
\end{enumerate}
We will need, in the proof of Theorem~\ref{embsmoothing}, $g_{{}_X} \geq 3$.
Since $g_{{}_X}=d+2g-1$, we will assume
\begin{enumerate}
\addtocounter{enumi}{2}
\item $p_a(\WY) \geq 3$.
\end{enumerate}
Observe that an obvious necessary condition for $\WY$ to be smoothable is $p_a(\WY) \geq 0$.
Thus, provided such a double cover $X$ there exists, the condition $p_a(\WY) \geq 3$ excludes only very few ribbons if $g=0$ or $g=1$, for then we must obviously have $d \geq 0$.\\
Moreover the existence of $X$ is a weak condition on $\SE$.
For instance the conditions (1) and (3) are verified if $d \geq \max\{g, -2g+4\}$ (we impose $d \geq g$ for $\SE^{-2}$ to be globally generated).
\begin{theorem}\label{embsmoothing}
Let $Y$ be a smooth irreducible projective curve and let $\SO_Y(1)$ be a very ample line bundle on $Y$.
Let $\SE$ be a line bundle on $Y$.
Assume that $\SO_Y(1)$ and $\SE \otimes \SO_Y(1)$ are nonspecial.
Let $\, \PP^s \subset \PP^r$ denote, respectively, the projective spaces of (one quotients) of $\,H^0(\SO_Y(1))$ and $H^0(\SO_Y(1)) \oplus H^0(\SE \otimes \SO_Y(1))$.
Assume that $r \geq 3$.
Let $Y \subset \PP^r$ be the embedding defined as the composition of the embedding $Y \subset \PP^s$ given by the complete linear series $H^0(\SO_Y(1))$ and $\PP^s \subset \PP^r$.\\
Assume that $\WY \subset \PP^r$ is a nondegenerate embedded ribbon over $Y \subset \PP^r$ with conormal bundle $\SE$.
Assume that $p_a(\WY) \geq 3$.
Assume that there is a smooth irreducible double cover $X \overset{\pi} \to Y$ with $\pi_*\SO_X/\SO_Y = \SE$.
Let $X \overset{\varphi}\to \PP^r$ be the morphism obtained as the composition of $\pi$ with the inclusion of $Y$ in $\PP^r$.
Then
\begin{enumerate}
\item there exists a smooth irreducible family $\SX$ proper and flat over a smooth pointed affine curve $(T, 0)$ and a $T$--morphism $\SX \overset{\Phi}\to \PP_T^r$ with the following properties:
\begin{enumerate}
\item the general fiber $\SX_t \overset{\Phi_t}\to \PP^r, t \neq 0,$ is a closed immersion of a smooth irreducible projective curve $\SX_t$,
\item the central fiber $\SX_0 \overset{\Phi_0}\to \PP^r$ is $X \overset{\varphi} \to \PP^r$; and
\end{enumerate}
\item the image of $\SX \overset{\Phi}\to \PP_T^r$ is a closed integral subscheme $\SY \subset \PP^r_T$ flat over $T$ with the following properties:
\begin{enumerate}
\item the ge\-neral fiber $\, \SY_t, t \neq 0,$ is a smooth irreducible projective nondegenerate curve with nonspecial hyperplane section in $\PP^r$,
\item the central fiber $\, \SY_0$ is $\,\WY \subset \PP^r$.
\end{enumerate}
\end{enumerate}
\end{theorem}
\begin{remark}
{\rm If $H^0(\SE \otimes \SO_Y(1))=0$ then $Y \subset \PP^r$ is nondegenerate.}{\hfill $\square$} 
\end{remark}
\begin{Proof}
(of Theorem~\ref{embsmoothing}) Let us denote $d =- \mathrm{deg} \, \SE$ and $g$ the genus of $Y$.
By assumption there is a smooth irreducible double cover $X \overset{\pi}\to Y$ with $\pi_*\SO_X=\SO_Y \oplus \SE$.
The curve $X$ is projective with genus $g_{{}_X}=d + 2g - 1$.
Thus $g_{{}_X} = p_a(\WY)$ and by assumption $g_{{}_X} \geq 3$.
On $X$ we consider the line bundle $L=\pi^*\SO_Y(1)$.
The hypothesis that $\SO_Y(1)$ and $\SE \otimes \SO_Y(1)$ are nonspecial implies that $L$ is nonspecial.\\
The natural map $H^0(\SO_Y(1)) \overset{\pi^*}\hookrightarrow H^0(L)$ admits a retraction $H^0(L)\overset{p}\to H^0(\SO_Y(1))$ obtained from the trace map $\pi_*\SO_X \to \SO_Y$ twisting by $\SO_Y(1)$ and taking global sections.
So we have
\begin{equation*}
H^0(L)= H^0(\SO_Y(1))\oplus H^0(\SE \otimes \SO_Y(1)).
\end{equation*}
We also see that the pullback of the surjective evaluation map $H(\SO_Y(1))\otimes \SO_Y \twoheadrightarrow \SO_Y(1)$ is the composition of $H^0(\SO_Y(1)) \otimes \SO_X \overset{\pi^* \otimes \,\mathrm{id}}\longrightarrow H^0(L) \otimes \SO_X$ and the evaluation map $H^0(L) \otimes \SO_X \to L$.\\
Now the composition of the surjective map $H^0(L)\otimes \SO_X \overset{p \otimes \, \mathrm{id}} \longrightarrow H^0(\SO_Y(1)) \otimes \SO_X$ and the pullback of the surjective evaluation map $H(\SO_Y(1))\otimes \SO_Y \twoheadrightarrow \SO_Y(1)$ defines a surjection $H^0(L) \otimes \SO_X \twoheadrightarrow L$ that gives a morphism $X \overset{\varphi}\to \PP^r$.
This morphism $\varphi$ is the composition $X \overset{\pi}\to Y \hookrightarrow \PP^s \hookrightarrow \PP^r$.\\
From Theorem~\ref{key.curves.0} there exists a first order infinitesimal deformation $\WX \overset{\wphi}\to \PP_{\Delta}^r$ of $\varphi$ such that the central fiber of the image of $\wphi$ is equal to the ribbon $\WY$. Let us denote $\WL=\wphi^*\SO_{\PP_{\Delta}^r}(1)$. Then $\WL$ restricts to $L$ on $X$.\\
The hypothesis that $\WY$ is nondegenerate in $\PP^r$ implies that if the image of $\wphi$ is contained in a closed subscheme given in $\PP_{\Delta}^r$ by a linear form with coeficients in $k[\epsilon]$, then the central fiber of this subscheme is $\PP^r$.
Therefore we see that, given a set of coordinates of $\PP^r$, the morphism $\WX \overset{\wphi}\to \PP_{\Delta}^r$ corresponds to the choice of a set of $r+1$ sections $\{\wl_0,\ldots,\wl_r\}$ in $\Gamma(\WL)$ such that generate $\WL$, all whose possible relations over $k[\epsilon]$ have nonunit coefficients and whose restriction to $H^0(L)$ is a set $\{l_0,\ldots,l_r\}$ of $r+1$ sections that generate $L$ and such that exactly $s+1$ of them are independent.
To this last set corresponds $X \overset{\varphi}\to \PP^r$.\\
We consider $\womega=\omega_{\WX/\Delta}$ and $\WL'=\WL \otimes \womega^{\otimes n}$, where $n$ is large enough so that $L'= L \otimes \omega_X^{\otimes n}$ is very ample, nonspecial and the complete linear series of $L'$ defines an embedding $X \hookrightarrow \PP^{r'}$ that determines a smooth point $[X']$ in the corresponding Hilbert scheme.
Let $H$ be the open, smooth and irreducible subset of this Hilbert scheme formed by smooth irreducible nondegenerate curves $C \subset \PP^{r'}$ of degree $d'=2 \, \mathrm{deg}\,\SO_Y(1)+n(2g_{{}_X}-2)$ and genus $g_{{}_X}=d+2g-1$.
Then $[X'] \in H$.
Since $n>>0$, then for every such curve $\SO_C(1)$ is nonspecial, the embedding of $C$ in $\PP^{r'}$ is defined by a complete series and defines a smooth point in its Hilbert scheme.
Moreover, since $L'$ is very ample and $H^1(L')=0$, also $\WL'$ is very ample relative to $\Delta$ and the embedding $X \hookrightarrow \PP^{r'}$ extends to an embedding $\WX \hookrightarrow \PP_{\Delta}^{r'}$.
So the image $\WX'$ of $\WX \hookrightarrow \PP_{\Delta}^{r'}$ is a flat family over $\Delta$ that corresponds to a tangent vector to $H$ at the Hilbert point $[X']$ of $X'$.
We can take the embedding $\WX \hookrightarrow \PP_{\Delta}^{r'}$ so that this tangent vector be nonzero.
Now, since $[X']$ is a smooth point in $H$, we can take a smooth irreducible affine curve $T$ in $H$ passing through $[X']$ with tangent direction the given tangent vector.\\
We can take this curve in such a way that all its points, except perhaps $[X']$, be placed in the open subset $U$ of $H$ constructed in the following way: $H$ admits a surjective morphism over $\SP_{d',g_{{}_X}}$, the coarse moduli of pairs consisting of a curve of genus $g_{{}_X}$ and a line bundle of degree $d'$ on the curve.
Denote $d_1=2 \,\mathrm{deg}\,\SO_Y(1)$ and consider also $\SP_{{d_1},g_{{}_X}}$ fibered over the fine part of the moduli $\SM_{g_{{}_X}}^0$.
Let $\SC^{(d_1)}$ be the scheme that represents the functor of relative effective Cartier divisors of relative degree $d_1$ over the universal curve $\SC_{g_{{}_X}}^0 \to \SM_{g_{{}_X}}^0$ (see \cite[4.1]{SemBourbaki232}).
By hypothesis $d_1 - g_{{}_X}=r \geq 3$, in particular $d_1 \geq g_{{}_X}$.
Therefore the morphism $\SC^{(d_1)} \to \SP_{{d_1},g_{{}_X}}$ is surjective.
Denote $\SC = \SC_{g_{{}_X}}^0 \times_{{}_{\SM_{g_{{}_X}}^0}} \SC^{(d_1)}$. Over $\SC$ there is a universal effective relative Cartier divisor $\SD$. Consider the line bundle $\SO_{\SC}(\SD)$ and let $\SC \overset{q}\to \SC^{(d_1)}$ be the (proper and flat) projection.
Then, by the theorem of base change and cohomology, at a point $(C, D)$ of $\SC^{(d_1)}$, consisting of a curve $C$ of genus $g_{{}_X}$ and a nonspecial divisor $D$ of degree $d_1$ on $C$, the fiber of the coherent sheaf $R^1q_*\SO_{\SC}(\SD)$ is isomorphic to $H^1(C, D)$ and the same is true near $(C, D)$. So there is a non--empty open set $W_1$ on $\SC^{(d_1)}$ formed by pairs consisting of a curve and a divisor whose associated line bundle is nonspecial.
Furthermore, if we restrict to $W_1$, then the support of the cokernel of the natural map $q^*q_*(\SO_{\SC}(\SD)) \to \SO_{\SC}(\SD)$ is outside of the inverse image of an open set $W_2 \subset W_1$. So we obtain an open set $W_2$ on $\SC^{(d_1)}$ formed by pairs consisting of a curve of genus $g_{{}_X}$ and an effective divisor of degree $d_1$ whose associated line bundle is nonspecial and globally generated.
In~\cite[5.1]{EisenbudHarris83} it is proved that if $r \geq 3$ then on a general smooth curve the general linear series of dimension $r$ has no base points and its associated map to $\PP^r$ is a closed immersion.
Moreover, $\SC^{(d_1)}$ is irreducible so $W_1$ dominates $\SM_{g_{{}_X}}^0$.
Therefore, the set $W_2$ is non--empty, and shrinking $W_2$ so that $q_*(\SO_{\SC}(\SD))$ is free of rank $r+1$ on $W_2$, we have a $W_2$--morphism 
$\, \SC_{g_{{}_X}}^0 \times_{{}_{\SM_{g_{{}_X}}^0}} W_2  \to \PP_{W_2}^r$.
Now (see e.g.~\cite[4.6.7]{EGA3-1}) the set $W$ formed by the points of $W_2$ such that the induced morphism on the fiber over the point is a closed immersion is open in $W_2$.
So we obtain an open set $W$ in $\SC^{(d_1)}$ formed by pairs consisting of a curve of genus $g_{{}_X}$ and an effective divisor of degree $d_1$ whose associated line bundle is nonspecial and very ample.
Also by~\cite[5.1]{EisenbudHarris83}, if we assume that $d_1 - g_{{}_X}=r \geq 3$ then the open set $W$ is non--empty.
Now, since $\SC^{(d_1)}$ is irreducible and $\SC^{(d_1)} \to \SP_{{d_1},g_{{}_X}}$ is surjective, we also obtain a non--empty open set in $\SP_{d_1,g_{{}_X}}$ formed by pairs consisting of a curve and a very ample nonspecial line bundle with as many global sections as $L$.
Moreover, twisting by $\omega^{\otimes n}$ we have an isomorphism between $\SP_{d_1,g_{{}_X}}$ and $\SP_{d',g_{{}_X}}$.
So we take the open set $\,U \subset H$ inverse image of the considered open set in $\SP_{d',g_{{}_X}}$.
Now, we take our curve $T$ with general point in this open set $U$.\\
Let $0 \in T$ denote the point corresponding to $X'$.
Over the pointed affine curve $(T, 0)$ we have a proper flat polarized family $(\SX,\SL')$ containing $(X, L')$ and $(\WX,\WL')$ as fibers over the point $0$ and the tangent vector to $T$ at $0$.
Now, twisting by the $-n$ power of the relative dualizing sheaf of the family we obtain a family $(\SX,\SL)$ proper and flat over $T$ whose central fiber is $(X,L)$, whose restriction to the tangent vector to $T$ at $0$ is $(\WX,\WL)$ and whose general member $(\SX_t,\SL_t)$ consists of a smooth irreducible projective curve of genus $g_{{}_X}$ and a very ample nonspecial line bundle $\SL_t$ with as many global sections as $L$ and degree $d_1=\,\mathrm{deg}\,L$.
Then we have $h^0(\SL_t)=r+1$ for every $t$ and we will show that, after shrinking $T$ if it is necessary, $\SL$ defines a $T$--morphism $\SX \to \PP_T^r$ whose fiber over the tangent vector to $T$ at $0$ is the initial morphism $\WX \overset{\wphi}\to \PP_{\Delta}^r$ and whose general fiber $\SX_t \overset{\varphi_t} \to \PP^r$ for $t \neq 0$ is a closed immersion given by the complete linear series of $\SL_t$.
Indeed, recall that the morphism $\WX \overset{\wphi}\to \PP_{\Delta}^r$ is associated with a surjective map $\SO_{\WX}^{r+1} \twoheadrightarrow \WL$ given by $r+1$ global sections $\{\wl_0,\ldots,\wl_r\}$ all whose possible relations over $k[\epsilon]$ have nonunit coefficients and whose restriction to $X$ is a set $\{l_0,\ldots,l_r\}$ of $r+1$ global sections of $L$ such that exactly $s+1$ of them are independent.
The restricted surjection $\SO_X^{r+1}\twoheadrightarrow L$ given by $\{l_0,\ldots,l_r\}$ defines the initial morphism $X \overset{\varphi}\to \PP^r$.
Now, we will obtain a $T$--morphism $\SX \to \PP_T^r$ extension of $\WX \overset{\wphi}\to \PP_{\Delta}^r$ if we can lift $\{\wl_0,\ldots,\wl_r\}$ to global sections $\{\mathsf{m}_0, \ldots,\mathsf{m}_r\}$ of $\SL$ such that the associated map $\SO_{\SX}^{r+1} \to \SL$ is surjective.\\ 
Let $\SX \overset{p}\to T$ be the (proper and flat) structural morphism.
The facts that $p$ is proper, $\SL$ is flat over $T$ and $H^1(\SX_t,\SL_t)=0$ for every $t \in T$ imply that $p_* \SL$ is a locally free sheaf of rank $r+1$ on $T=\mathrm{Spec}\, R$ and ``the formation of $p_*$ commutes with base extension'' so we have $\Gamma(\SL) \otimes_R k[\epsilon]/\epsilon k[\epsilon] = \Gamma(L)$ and $\Gamma(\SL) \otimes_R k[\epsilon] = \Gamma(\WL)$.
After shrinking $T$, we can assume that $\mathsf{M}=\Gamma(\SL)$ is a free $R$--module of rank $r+1$.\\
We prove that the map $\mathsf{M} \to \Gamma(\WL)$ is surjective.
Twisting by $\SL$ the short exact sequence associated with the inclusion $\WX \subset \SX$ and pushing--down to $T$, we obtain an exact sequence
\begin{equation*}
\xymatrix@C-5pt{
0 \ar[r] & p_*\SL(-2X) \ar[r] & p_* \SL \ar[r] & \wpp_* \WL \ar[r]  & R^1p_*\SL(-2X),}
\end{equation*}
where $\WX \overset{\wpp}\to \Delta$ is the structural morphism.
Now, shrinking $T$, we can assume that $\SO_{\SX}(-X)$ is isomorphic to $\SO_{\SX}$ and thus we see that $R^1p_*\SL(-2X)$ vanishes from the fact that $\SL$ induces nonspecial line bundles on every fiber.
So we can lift $\{\wl_0,\ldots,\wl_r\}$ to sections $\{\mathsf{m}_0,\ldots, \mathsf{m}_r \}$.
The sections $\{\mathsf{m}_0,\ldots, \mathsf{m}_r \}$ define a map $\SO_{\SX}^{r+1} \to \SL$ whose cokernel vanishes at $0$.
Therefore, shrinking $T$, we can assume that $\SO_{\SX}^{r+1} \to \SL$ is surjective.
Thus we have obtained a surjection $\SO_{\SX}^{r+1} \to \SL$ that defines a $T$--morphism $\SX \overset{\Phi}\to \PP_T^r$ whose $\Delta$--fiber is $\WX \overset{\wphi}\to \PP_{\Delta}^r$.\\
The section $\mathsf{m}_0 \wedge \cdots \wedge \mathsf{m}_r $ of $\wedge^{r+1} \mathsf{M}$ corresponds, after a choice of basis in $\mathsf{M}$, to an element $\mathsf{d} \in R$.
I claim that $\mathsf{d}\neq 0$. Indeed, we see this by showing that $\mathsf{d}$ does not vanish at order $n=r-s$.
If $n=0$ then $\{l_0,\ldots,l_r\}$ are independent so $\mathsf{d}$ does not vanish at $0 \in T$.
Assume $n\geq 1$.
The $k[\epsilon]$--module $\Gamma(\WL)=\mathsf{M} \otimes_R k[\epsilon]$ is free so we have $\Gamma(\WL)=\Gamma(L)\oplus \Gamma(L)\, \epsilon$.
Therefore we can write $\wl_i =l_i + m_i \epsilon$ where $m_i \in \Gamma(L)$.
The vanishing of $\mathsf{d}$ at order $n$ is equivalent to
\begin{equation}\label{bigsum}
\sum_{0\leq i_1 < \cdots < i_n \leq r} l_0 \wedge \cdots \wedge l_{i_1-1}\wedge m_{i_1} \wedge l_{i_1+1} \wedge \cdots \wedge l_{i_n-1}\wedge m_{i_n}\wedge l_{i_n+1}\wedge \cdots \wedge l_r \, =0.  
\end{equation}
From~\eqref{bigsum} we obtain a $k[\epsilon]$--linear relation among the sections $\{l_0+m_0 \epsilon,\ldots,l_r +m_r \epsilon\}$ such that some of its coefficients is a unit in $k[\epsilon]$.
The existence of this linear relation implies that the central fiber $\,(\mathrm{im}\,\widetilde{\varphi})_0$ is degenerate.
This is contrary to our hypothesis that the ribbon $\WY$ is nondegenerate.
So the equality~\eqref{bigsum} does not happen and therefore $\mathsf{d}$ does not vanish at order $n=r-s \,$ as we wanted to show.\\
Therefore, shrinking $T$, we can assume that the $r+1$ elements $\{\mathsf{m}_0,\ldots, \mathsf{m}_r \}$ of $\mathsf{M}$ induce a basis in $H^0(\SL_t)$ for every $0 \neq t \in T$.\\
Thus we obtain for every $0 \neq t \in T$ a surjection $\SO_{\SX_t}^{r+1} \twoheadrightarrow \SL_t$ given by a basis of $H^0(\SL_t)$ and for $t=0$ the surjection $\SO_X^{r+1}\twoheadrightarrow L$ whose associated morphism is $X\overset{\varphi}\to \PP^r$.
This is a flat family of morphisms $\SX_t \overset{\Phi_t} \to \PP^r$ whose central fiber is $X\overset{\varphi}\to \PP^r$, whose general fiber is a closed immersion associated with a complete linear series and whose $\Delta$--fiber is $\WX \overset{\wphi}\to \PP_{\Delta}^r$.\\
Let $\SY$ be the image of the $T$--morphism $\SX \overset{\Phi}\to \PP_T^r$.
The total family $\SX$ is smooth and irreducible so $\SY$ is integral.
Furthermore, $\Phi$ is a closed immersion over $T-0$ since $\Phi_t$ is a closed immersion for every $t \in T-0$ (see e.g. \cite[4.6.7]{EGA3-1}).
Therefore for $t \in T-0$ we have the equality $\SY_t = \,\mathrm{im}\,(\Phi_t)$.
Finally, the facts that $T$ is an integral smooth curve and $\SY$ is integral and dominates $T$ imply that $\SY$ is flat over $T$.
So the fiber $\SY_0$ of $\SY$ at $0 \in T$ is the flat limit of the images of $\SX_t \overset{\Phi_t} \to \PP^r$ for $t \neq 0$.
Moreover, this fiber $\SY_0$ contains the central fiber $\,(\mathrm{im}\,\widetilde{\varphi})_0$ of the image of $\wphi$ and since both the fiber $\SY_0$ and the fiber $(\mathrm{im}\,\widetilde{\varphi})_0$ have the same degree and the same arithmetic genus they are equal.
\end{Proof}\\
As consequence of Theorem~\ref{embedding} and Theorem~\ref{embsmoothing} we obtain the smoothing of ribbons of arithmetic genus greater than or equal to $3$.
\begin{theorem}\label{abssmoothing}
Let $Y$ be a smooth irreducible projective curve and let $\SE$ be a line bundle on $Y$.
Assume that there is a smooth irreducible double cover $X \overset{\pi} \to Y$ with $\pi_*\SO_X/\SO_Y = \SE$.
Then every ribbon $\WY$ over $Y$ with conormal bundle $\SE$ and arithmetic genus $p_a(\WY) \geq 3$ is smoothable.{\hfill $\square$}
\end{theorem}
The conditions (1) there is a non--zero effective reduced divisor on $Y$ with associated line bundle $\SE^{-2}$ and (3) $p_a(\WY) \geq 3$ are verified if $d \geq \max \{g, -2g+4\}$ so we obtain:
\begin{cor}\label{abssmoothing2}
Let $Y$ be a smooth irreducible projective curve of genus $g$.
\begin{enumerate}
\item Let $\SE$ be a line bundle on $Y$ and $d = -\mathrm{deg} \, \SE$.
If $d \geq \max \{g, -2g+4\}$ then every ribbon over $Y$ with conormal bundle $\SE$ is smoothable.
\item Assume that $g \geq 2$.
Let $\SE$ be a line bundle on $Y$ such that $\SE^{-2}=\SO_Y$ and $H^0(\SE)=0$.
Then every ribbon over $Y$ with conormal bundle $\SE$ is smoothable.{\hfill $\square$}
\end{enumerate}
\end{cor}
For ribbons over an elliptic curve or ribbons over $\PP^1$, we can also apply Theorem~\ref{embsmoothing} for ribbons embedded like in Remark~\ref{elliptic-rational}.
In these cases we obtain the following embedded smoothing results for ribbons supported over a nondegenerate embedding of its reduced part.
\begin{cor}\label{elliptic-smoothing}
Let $Y \subset \PP^{d-1}$ be an elliptic normal curve of degree $d \geq 5$.
Let $\SE$ be a line bundle of degree $-d$ such that $\SE^{-1}$ is not isomorphic to $\SO_Y(1)$.
Then for every ribbon $\WY$ over $Y$ with conormal bundle $\SE$ embedded in $\PP^{d-1}$ with support on $Y \subset \PP^{d-1}$ there exists a closed integral subscheme $\SY \subset \PP^{d-1} \times T$ flat over a smooth pointed affine curve $T$ whose ge\-neral fiber is a smooth irreducible projective nondegenerate curve of genus $d+1$ with nonspecial hyperplane section in $\PP^{d-1}$ and whose central fiber is $\WY \subset \PP^{d-1}$.
Moreover, in this conditions, every ribbon $\WY$ over $Y$ with conormal bundle $\SE$ admits an embedding in $\PP^{d-1}$ with support on $Y \subset \PP^{d-1}$.{\hfill $\square$}
\end{cor}
\begin{cor}\label{rational-smoothing}
Let $Y \subset \PP^{d-1}$ be a rational normal curve of degree $d-1$.        
Assume that $d \geq 4$.
Then for every ribbon $\WY$ over $\PP^1$ with conormal bundle $\SO_{\PP^1}(-d)$ embedded in $\PP^{d-1}$ with support on $Y \subset \PP^{d-1}$ there exists a closed integral subscheme $\SY \subset \PP^{d-1} \times T$ flat over a smooth pointed affine curve $T$ whose ge\-neral fiber is a smooth irreducible projective nondegenerate curve of genus $d-1$ with nonspecial hyperplane section in $\PP^{d-1}$ and whose central fiber is $\WY \subset \PP^{d-1}$.
Moreover, in this conditions, every ribbon $\WY$ over $\PP^1$ with conormal bundle $\SO_{\PP^1}(-d)$ admits an embedding in $\PP^{d-1}$ with support on $Y \subset \PP^{d-1}$.{\hfill $\square$}
\end{cor}
\bibliographystyle{amsalpha}

\providecommand{\bysame}{\leavevmode\hbox to3em{\hrulefill}\thinspace}
\providecommand{\MR}{\relax\ifhmode\unskip\space\fi MR }
\providecommand{\MRhref}[2]{%
  \href{http://www.ams.org/mathscinet-getitem?mr=#1}{#2}
}
\providecommand{\href}[2]{#2}

\end{document}